\documentclass[preprint,authoryear,3p,times]{elsarticle}

\usepackage{amscd,amssymb,amsmath,amsthm}
\usepackage{graphicx,verbatim,lscape,rotating}
\usepackage{color,float,caption,subcaption}
\usepackage{booktabs,tabularx,multirow,array}
\usepackage{dsfont}
\usepackage{here}
\usepackage{float}

\newcommand{\bfzero}{\mathbf {0}}

\newcommand{\bfb}{\mathbf {b}}

\newcommand{\bfe}{\mathbf {e}}

\newcommand{\bfr}{\mathbf {r}}

\newcommand{\bfz}{\mathbf {z}}
\newcommand{\bfA}{\mathbf {A}}
\newcommand{\bfB}{\mathbf {B}}
\newcommand{\bfC}{\mathbf {C}}

\newcommand{\bfE}{\mathbf {E}}

\newcommand{\bfI}{\mathbf {I}}

\newcommand{\bfM}{\mathbf {M}}

\newcommand{\bfR}{\mathbf {R}}

\newcommand{\bfU}{\mathbf {U}}

\newcommand{\bfW}{\mathbf {W}}
\newcommand{\bfX}{\mathbf {X}}
\newcommand{\bfY}{\mathbf {Y}}
\newcommand{\bfZ}{\mathbf {Z}}
\newcommand{\cd}{\stackrel{d}{\to}}
\newcommand{\cp}{\stackrel{p}{\to}}
\newcommand{\cas}{\stackrel{a.s.} {\to}}

\newcommand{\bfbeta}{\mbox{\boldmath$\beta$}}

\newcommand{\bfepsilon}{\mbox{\boldmath$\epsilon$}}

\newcommand{\bfmu}{\mbox{\boldmath$\mu$}}
\newcommand{\bfeta}{\mbox{\boldmath$\eta$}}

\newcommand{\bfxi}{\mbox{\boldmath$\xi$}}
\newcommand{\bfSigma}{\mbox{\boldmath$\Sigma$}}

\newcommand{\bfOmega}{\mathbf{\Omega}}
\newcommand{\bfPhi}{\mathbf{\Phi}}
\newcommand{\bfPsi}{\mathbf{\Psi}}
\newcommand{\bfTheta}{\mathbf{\Theta}}

\newcommand{\tr}{\mbox{tr}}

\newcommand{\var}{\mbox{var}}
\newcommand{\cov}{\mbox{cov}}

\newtheorem{theorem}{Theorem}

\newtheorem{remark}{Remark}
\numberwithin{theorem}{section}
\numberwithin{corollary}{section}
\numberwithin{lemma}{section}
\numberwithin{remark}{section}

\def\vec{\mbox{\textnormal{vec}}}

\usepackage{lineno}
\modulolinenumbers[5]

\begin{document}

\begin{frontmatter}

%% Title, authors and addresses

%% use the tnoteref command within \title for footnotes;
%% use the tnotetext command for the associated footnote;
%% use the fnref command within \author or \address for footnotes;
%% use the fntext command for the associated footnote;
%% use the corref command within \author for corresponding author footnotes;
%% use the cortext command for the associated footnote;
%% use the ead command for the email address,
%% and the form \ead[url] for the home page:
%%
%% \title{Title\tnoteref{label1}}
%% \tnotetext[label1]{}
%% \author{Name\corref{cor1}\fnref{label2}}
%% \ead{email address}
%% \ead[url]{home page}
%% \fntext[label2]{}
%% \cortext[cor1]{}
%% \address{Address\fnref{label3}}
%% \fntext[label3]{}

\title{Estimating weak periodic vector autoregressive time series}

%% use optional labels to link authors explicitly to addresses:
%% \author[label1,label2]{<author name>}
%% \address[label1]{<address>}
%% \address[label2]{<address>}
\author[franche]{Yacouba Boubacar Ma\"{i}nassara \corref{cor1}} \ead{yacouba.boubacar\_mainassara@univ-fcomte.fr}
\cortext[cor1]{Corresponding author}
\address[franche]{ Universit\'e Bourgogne Franche-Comt\'e, Laboratoire de math\'ematiques de Besançon, UMR CNRS 6623, 16 Route de Gray, 25030 Besançon, France}
\author[BSE,ecreb]{Eugen Ursu}\ead{eugen.ursu@u-bordeaux.fr}
\address[BSE]{Universit\'e de Bordeaux, BSE - UMR CNRS 6060, 16 Avenue L\'eon Duguit, B\^{a}t. H2, 33608 Pessac CEDEX, France}
\address[ecreb]{ECREB, Faculty of Economics and Business Administration, West University of Timisoara, 16 J.H. Pestalozzi Street, 300115, Timisoara, Romania}
%\address{}

\begin{abstract}
This article develops the asymptotic distribution of the least squares estimator of the model parameters in periodic vector autoregressive time series models (hereafter PVAR) with uncorrelated but dependent innovations. When the innovations are dependent, this asymptotic distributions can be quite different from that of PVAR models with independent and identically distributed (iid for short) innovations developed in~\citet{UD09}. Modified versions of the Wald tests are proposed for testing linear restrictions on the parameters. These asymptotic results are illustrated by Monte Carlo experiments. An application to a bivariate real financial data is also proposed.
\end{abstract}


\end{frontmatter}

%%
%% Start line numbering here if you want
%%
% \linenumbers

%% main text
\section{Introduction}
\noindent Many phenomena observed over time are subject to seasonal effects, which
are variations occurring at specific regular time intervals. The autoregressive
integrated moving average (ARIMA) model could be modified by employing the seasonal
differencing operator: if considered period magnitude is $s$, this operator subtracts from each observation the corresponding value at $s$ previous time instants. The result is the seasonal autoregressive integrated moving average (SARIMA) model developed originally by~\citet{BJ70}. This way of proceeding has been proven useful when mean for a given season is not stationary across years~\citep{HM94}. However, it turns out that many seasonal time series cannot be filtered to achieve second-order stationarity due to the correlation structure of these time series with the season~\citep{Ve85b}. For this reason a different procedure of accounting for seasonality has been proposed in literature, leading to periodic models. The use of periodic models appears to be well-suited to deal with many real life phenomena characterized by a seasonal behavior: climatology~\citep{LLL10}, hydrology~\citep{Ve85a}, macroeconomics~\citep{FP04} and engineering~\citep{SDS13}. Multivariate models are expected to be more useful in practice, since most real-life situations involve several variables and vector time series. The maximum likelihood estimation~\citep{Lu05} and the least squares (LS) method~\citep{UD09} are efficient methods to estimate the PVAR models. However, the innovations in these PVAR models have the iid property. We refer to these as strong PVAR models, by opposition to weak PVAR models where the innovations are only uncorrelated.

In recent years, a large part of the time series and econometric literature was devoted to weaken the independence innovation assumption. This independence assumption is restrictive because it excludes conditional heteroscedasticity and other forms of nonlinearity~\citep{frz05}. Another argument in favor for considering the weak PVAR models comes from~\citet{WGVM05} as they found evidences of the existence of autoregressive conditional heteroskedastic (ARCH) effects in modelling daily streamflow series in China. They argued that, the strong periodic autoregressive (hereafter PAR) models would perform better than the SARIMA model for capturing the ARCH effect in monthly flow series, but insufficient to fully capture the ARCH effects in daily flow series. \citet{frs11}~investigate the asymptotic properties of weighted least squares (WLS) estimation for causal and invertible periodic autoregressive moving average (PARMA) models with uncorrelated but dependent errors.~\citet{FR07} proposed a method to adjust the critical values of the portmanteau test for multiple autoregressive time series models with nonindependent innovations.

This article is organized as follows. In Section~\ref{model}, the weak PVAR model is introduced and the asymptotic properties of the least squares estimators are given in Section~\ref{result1}. An example of analytic computation of the asymptotic variance matrices is also given  in Section~\ref{result1}. Two weakly consistent estimators of the asymptotic variance matrix are proposed in Section~\ref{result2}.
In Section \ref{test}  it is shown how the standard Wald test must be adapted in the weak PVAR case in order to test for general linearity constraints. This section is also of interest in the univariate framework because, to our knowledge, this test has not been studied for weak PAR models. In Section~\ref{simul}, some simulation results are reported, and in Section~\ref{real}, an application to the daily returns of two European stock market indices, CAC 40 (Paris) and
DAX (Frankfurt), is made. Finally, Section~\ref{conclusion} offers some concluding remarks.

%%% ----------------------------------------------------------------------
\section{Weak periodic vector autoregressive time series models}\label{model}
\noindent In this section, we present principal results on the least squares estimators in the unconstrained and constrained case of the weak PVAR model.

Let $\bfY = \{ \bfY_t, t \in \mathbb{Z} \}$ be a stochastic process, where $$\bfY_t = (Y_t(1),\ldots,Y_t(d))^{\top}$$ represents a random vector of dimension $d$.
The process $\bfY$ is a PVAR process of order $p(\nu)$, $\nu \in \{1,\ldots,s\}$ (s is a predetermined value), if there exist
$d\times d$ matrices $\bfPhi_k(\nu) = \left(  \Phi_{k,ij}(\nu) \right)_{i,j=1,\ldots,d}$, $k=1,\ldots,p(\nu)$ such that
\begin{equation}
\label{pvar}
 \bfY_{ns+\nu} = \sum_{k=1}^{p(\nu)} \bfPhi_k(\nu) \bfY_{ns+\nu-k} + \bfepsilon_{ns+\nu}.
\end{equation}

The process $\bfepsilon:=(\bfepsilon_{t})_t=(\bfepsilon_{ns+\nu})_{n\in \mathbb Z}$ can be interpreted as in \cite{frs11} as the linear innovation of $\bfY:=(\bfY_t)_t=(\bfY_{ns+\nu})_{n\in \mathbb Z}$, \textit{i.e.} $\bfepsilon_t=\bfY_t-\mathbb{E}[\bfY_t|\mathcal{H}_\bfY(t-1)]$, where $\mathcal{H}_\bfY(t-1)$ is the Hilbert space generated by $(\bfY_u, u<t)$. The innovation process $\bfepsilon$ is assumed to be a stationary sequence satisfying
\begin{itemize}
\item[\hspace*{1em} {\bf (A0):}]
\hspace*{1em} $\mathbb{E}\left[ \bfepsilon_t\right]=0, \ \mathrm{Var}\left(\bfepsilon_t\right)=\Sigma_{\bfepsilon}(\nu) \text{ and } \mathrm{Cov}\left(\bfepsilon_t,\bfepsilon_{t-h}\right)=0$ for all $t\in\mathbb{Z}$ and all $h\neq 0$. The covariance matrix $\Sigma_{\bfepsilon}(\nu)$ is assumed to be non-singular.
\end{itemize}
Under the above assumptions the process $(\bfepsilon_{ns+\nu})_{n\in\mathbb{Z}}$ is called a weak multivariate periodic white noise.
An example of weak multivariate periodic white noise is the multivariate periodic  generalized autoregressive conditional heteroscedastic (MPGARCH) model (see for instance  \cite{Bi18}).
Many others univarites examples can also be find in \cite{frs11} and can be extended to multivariate  periodic white noise.

It is customary to say that  $(\bfY_{ns+\nu})_{n\in\mathbb{Z}}$ is a strong PVAR representation and we will do this henceforth if in  \eqref{pvar} $(\bfepsilon_{ns+\nu})_{n\in\mathbb{Z}}$
is a  strong multivariate periodic white noise, namely an iid sequence of random variables with mean $0$ and common variance matrix. A strong white noise is obviously  a weak white noise because independence entails uncorrelatedness. Of course the converse is not true.  In contrast with this previous definition, the representation  \eqref{pvar} is called a weak PVAR if  no additional assumption is made on $(\bfepsilon_{ns+\nu})_{n\in\mathbb{Z}}$, that is if $(\bfepsilon_{ns+\nu})_{n\in\mathbb{Z}}$ is only a weak periodic white noise (not
necessarily iid).
It is clear from these definitions that the following inclusions hold:
$$\left\{\text{strong PVAR} \right\}\subset\left\{\text{weak PVAR}\right\}.$$
Nonlinear models are becoming more and more employed because numerous real time series exhibit nonlinear dynamics. For instance conditional heteroscedasticity can not be generated by PVAR models with iid noises.\footnote{To cite few univariates examples of nonlinear processes, let us mention the generalized autoregressive conditional heteroscedastic (GARCH), the self-exciting threshold autoregressive (SETAR), the smooth transition autoregressive (STAR), the exponential autoregressive (EXPAR), the bilinear, the random coefficient autoregressive (RCA), the
functional autoregressive (FAR) (see \cite{fz19}, \cite{Tong1990} and \cite{FY2008} for references on these nonlinear time series models).}
As mentioned by \cite{fz05,fz98} in the case of autoregressive moving average (ARMA) models, many important classes of nonlinear processes admit weak ARMA representations.

The main issue with nonlinear models is that they are generally hard to identify and implement. These technical difficulties certainly explain the reason why the asymptotic theory of  PVAR model estimation is mainly limited to the strong PVAR model.

To derive some basic properties, it is convenient to write the model~(\ref{pvar}) in VAR representation:
\begin{equation}
\label{VARrepre}
\bfPhi_0^{\ast} \bfY_n^{\ast} = \sum_{k=1}^{p^*} \bfPhi_k^{\ast} \bfY_{n-k}^{\ast} +
                                \bfepsilon_n^{\ast},
\end{equation}
where
$\bfY_n^{\ast}= ( \bfY_{ns+s}^{\top}, \bfY_{ns+s-1}^{\top}, \ldots, \bfY_{ns+1}^{\top} )^{\top}$
and
$\bfepsilon_n^{\ast}= ( \bfepsilon_{ns+s}^{\top}, \bfepsilon_{ns+s-1}^{\top}, \ldots, \bfepsilon_{ns+1}^{\top} )^{\top}$
are
$(ds) \times 1$
random vectors.
The autoregressive model order in~(\ref{VARrepre})
is given by
$p^{\ast} = \lfloor p/s \rfloor$,
where
$\lfloor x \rfloor$
denotes the smallest integer greater than or equal to the real number $x$.
The
matrix
$\bfPhi_0^{\ast}$,
and the
autoregressive coefficients
$\bfPhi_k^{\ast}$, $k=1,\ldots,p^{\ast}$,
all of dimension
$(ds) \times (ds)$,
are given by the non-singular matrix:
\[
\bfPhi_0^{\ast}=
  \left [ \begin{array}{cccccc}
\bfI_d & -\bfPhi_1(s) & -\bfPhi_2(s) & \ldots & -\bfPhi_{s-2}(s) & -\bfPhi_{s-1}(s)\\
\bfzero & \bfI_d & -\bfPhi_1(s-1) & \ldots & -\bfPhi_{s-3}(s-1) & -\bfPhi_{s-2}(s-1)\\
\vdots &  & & \ddots & & \vdots  \\
\bfzero & \bfzero & \bfzero & \ldots & \bfI_d & -\bfPhi_1(2)\\
\bfzero & \bfzero & \bfzero & \ldots & \bfzero & \bfI_d
          \end{array}
  \right ],
\]
where
$\bfI_d$
denotes the
$d \times d$
identity matrix,
and:
\[
\bfPhi_k^{\ast}=
   \left [ \begin{array}{cccc}
\bfPhi_{ks}(s) & \bfPhi_{ks+1}(s) & \ldots & \bfPhi_{ks+s-1}(s) \\
\bfPhi_{ks-1}(s-1) & \bfPhi_{ks}(s-1) & \ldots & \bfPhi_{ks+s-2}(s-1)\\
\vdots & & \ddots & \vdots \\
\bfPhi_{ks-s+1}(1) & \bfPhi_{ks-s+2}(1) &  \ldots & \bfPhi_{ks}(1)
           \end{array}
    \right ],
\]
where
$k = 1,2,\ldots,p^{\ast}$
and
$\bfPhi_k(\nu) = \bfzero$, $k > p$.

\noindent From \eqref{VARrepre} we can in principle deduce the properties of weak PVAR parameters estimation, identification and validation from existing results on parameters estimation, identification and validation of the weak VARMA (Vector ARMA) models (see for instance \cite{BMF11,BM2012,BMK2016,BM11,BMS2018}).
Therefore we have preferred to work in the PVAR setting for various reasons. Firstly, in particular the results obtained directly in terms of the univariate PAR representation are more directly usable because fewer parameters are involved and their estimation is easier (see \cite{frs11} for more details).
Secondly,  the number of parameters in \eqref{VARrepre} is very huge, which entails statistical difficulties. Finally the VAR representation \eqref{VARrepre} is so-called structural form and is not standard when the matrix $\bfPhi_0^{\ast}\neq \bfI_{ds}$. The structural PVAR representation \eqref{VARrepre} can be rewritten in a standard reduced PVAR form if the matrix $\bfPhi_0^{\ast}$ is non singular, by multiplying \eqref{VARrepre} by $\bfPhi_0^{\ast-1}$ and introducing the innovation process $e_n=\bfPhi_0^{\ast-1}\bfepsilon_n^\ast$, with non singular variance $\bfPhi_0^{\ast-1}\text{var}(\bfepsilon_n^\ast) \left(\bfPhi_0^{\ast-1}\right)^\top$. This rescaling operation complicates the interpretation of the estimated parameters and the derivation of their statistical properties in the original scale, since the covariance matrix of the error term of the standard VARMA model now depends on the autoregressive parameters. The structural form \eqref{VARrepre} is mainly used in econometric to introduce instantaneous relationships between economic variables. The reduced form  is more practical from a statistical viewpoint, because it gives the forecasts of each component of $({\bfY}_n^\ast)$ according to the past values of the set of the components. The above discussion shows that the PVAR representation \eqref{VARrepre} is not unique, that is, a given process $({\bfY}_n^\ast)$ can be written in reduced form or in structural form by pre-multiplying by any non singular $(ds\times ds)$ matrix. Of course, in order to ensure the uniqueness of this  representation, constraints are necessary for the identifiability of the $(p^\ast+2)d^2s^2$ elements of the matrices involved in the PVAR equation \eqref{VARrepre}.

\noindent Let $\det (\bfA)$ be the determinant of the squared matrix $\bfA$.
Using general properties of VAR models,
it follows that the multivariate stochastic process
$\{ \bfY_t^{\ast} \}$
is causal if:
\begin{itemize}
\item[\hspace*{1em} {\bf (A1):}]
\hspace*{1em} $\det \left(\bfPhi_0^{\ast} - \bfPhi_1^{\ast} z - \ldots - \bfPhi_{p^{\ast}} z^{p^{\ast}} \right) \neq 0,$
for all complex numbers $z$ satisfying the condition
$|z|\leq 1$.
\end{itemize}
%\begin{equation}
%\label{causal}
%  \det \left(\bfPhi_0^{\ast} - \bfPhi_1^{\ast} z - \ldots - \bfPhi_{p^{\ast}} z^{p^{\ast}} \right) \neq 0,
%\end{equation}
%for all complex numbers $z$ satisfying the condition
%$|z|\leq 1$,
%where
%$\det (\bfA)$
%denotes the determinant of the squared matrix $\bfA$.
Under Assumption {\bf (A1)}, there exists a sequence of constant matrices $(\bfC_i(\nu))_{i\geq0}$ such that, for $\nu=1,2,\dots,s$,
$\sum_{i=0}^{\infty}\|\bfC_i(\nu)\|<\infty$ with $\bfC_0(\nu)=\bfI_d$ and
\begin{equation}
\label{MAinf}
\bfY_{ns+\nu}=\sum_{i=0}^{\infty}\bfC_i(\nu)\bfepsilon_{ns+\nu-i},
\end{equation}
where the sequence of matrices  $\|\bfC_i(\nu)\|\to 0$ at a geometric rate as $i\to\infty$. The $\|\bfA\|$ denotes the Euclidean norm of the matrix $\bfA$,
that is $\|\bfA\| = \{  \tr (\bfA \bfA^{\top}) \}^{1/2}$, with $\tr(\bfB)$ being the trace of the squared matrix $\bfB$.

To establish the consistency  of the least squares estimators, an additional assumption is needed.
\begin{itemize}
\item[\hspace*{1em} {\bf (A2):}]
\hspace*{1em}  The $ds$-dimensional process $\left(\boldsymbol{\epsilon}^{\ast}_n\right)_{n\in \mathbb{Z}}$ is ergodic and strictly  stationary.
\end{itemize}
Note that Assumption {\bf (A2)} is entailed by an iid assumption on  $\boldsymbol{\epsilon}^{\ast}_n$, but not by Assumption {\bf (A0)}.

For the asymptotic normality of  least squares estimators, additional assumptions are also required.
To control the serial dependence of the stationary process $(\boldsymbol{\epsilon}^{\ast}_n)_{n\in\mathbb{Z}}$, we introduce the strong mixing coefficients $\alpha_{\boldsymbol{\epsilon}^{\ast}}(h)$ defined by
$$\alpha_{\boldsymbol{\epsilon}^{\ast}}\left(h\right)=\sup_{A\in\mathcal F^n_{-\infty},B\in\mathcal F_{n+h}^{+\infty}}\left|\mathbb{P}\left(A\cap
B\right)-\mathbb{P}(A)\mathbb{P}(B)\right|,$$
where $\mathcal F_{-\infty}^n=\sigma (\boldsymbol{\epsilon}^{\ast}_u, u\leq n )$ and $\mathcal F_{n+h}^{+\infty}=\sigma (\boldsymbol{\epsilon}^{\ast}_u, u\geq n+h )$.

We use $\|\cdot\|$ to denote the Euclidean norm  of a vector.
We will make an integrability assumption on the moment of the noise and a summability condition on the strong mixing coefficients $(\alpha_{\boldsymbol{\epsilon}^{\ast}}(k))_{k\ge 0}$.
\begin{itemize}
\item[\hspace*{1em} {\bf (A3):}]
\hspace*{1em}
$\text{We have }\mathbb{E}\|{\boldsymbol{\epsilon}}^{\ast}_n\|^{4+2\kappa}<\infty\text{ and }
\sum_{k=0}^{\infty}\left\{\alpha_{{\boldsymbol{\epsilon}}^{\ast}}(k)\right\}^{\frac{\kappa}{2+\kappa}}<\infty\text{ for some } \kappa>0.$
\end{itemize}

%--------------------------------------------------------------------------------------------------------------------------------------
\section{Unconstrained least squares estimators and least squares estimation with linear constraint on the parameters.\\}\label{result1}
\noindent In this section, we study the asymptotic properties of least squares estimators from a causal PVAR model.
Let $\bfbeta(\nu) = (\vec^{\top}\{ \bfPhi_1(\nu) \},\ldots, \vec^{\top}\{ \bfPhi_{p(\nu)}(\nu) \} )^{\top}$ be a $\{ d^2 p(\nu) \} \times 1$ vector of parameters, where $\vec(\bfA)$ corresponds to the vector obtained by stacking the columns of $\bfA$~\citep[Chapter~16.3]{Ha97}
The PVAR model in~(\ref{pvar}) has $d^2 \sum_{\nu=1}^{s} p(\nu)$ autoregressive parameters $\bfPhi_k(\nu)$, $k=1,\ldots,p(\nu)$, $\nu = 1,\ldots,s$,
and $s$ additional $d \times d$ covariance matrices $\bfSigma_{\bfepsilon}(\nu)$,
$\nu = 1,\ldots,s$. For multivariate processes, the number of parameters can be quite large; for vector periodic processes, the inflation of parameters is due to the $s$ seasons. For example, in the case of bivariate monthly data where $d=2$, $s=12$, and, in the simplest case $p(\nu) \equiv 1$, this means that 48 independent autoregressive parameters must be estimated (by comparison, a traditional VAR(1) process relies on four independent parameters). In view of these considerations,
we consider estimation in the unrestricted case but also in the situation where the parameters of the same season $\nu$ satisfy the relation:
\begin{equation}
\label{constraints}
\bfbeta(\nu) = \bfR(\nu) \bfxi(\nu) + \bfb(\nu),
\end{equation}
where $\bfR(\nu)$ is a known
$\{ d^2 p(\nu) \} \times K(\nu)$
matrix of rank
$K(\nu)$, $\bfb(\nu)$ a known $\{ d^2 p(\nu) \} \times 1$
vector and
$\bfxi(\nu)$
represents a
$K(\nu) \times 1$
vector of unknown parameters.
Letting $\bfR(\nu) = \bfI_{d^2p(\nu)}$, $\bfb(\nu) = \bfzero$, $\nu = 1,\ldots,s$
give what we call the full unconstrained case.
In general, the matrices
$\bfR(\nu)$
and the vectors
$\bfb(\nu)$
allow for linear constraints on the parameters of the same season
$\nu$,
$\nu = 1,\ldots,s$.

This linear constraint includes the important special case of parameters set to zero
on certain components of $\bfPhi_k(\nu)$, $\nu = 1,\ldots,s$. In practice, a two-step procedure could consist of fitting a full unconstrained model, and, in a second stage of inference, the estimators which are statistically not significant could be considered known zero parameters, providing frequently more parsimonious models.

Consider the time series data $\bfY_{ns+\nu}$, $n=0,1,\ldots,N-1$, $\nu = 1,\ldots,s$, giving a sample size equal to $n = Ns$. Let
\begin{eqnarray}
\label{Znu}
\bfZ(\nu) &=& \left( \bfY_{\nu}, \bfY_{s+\nu}, \ldots, \bfY_{(N-1)s+\nu} \right), \\
\label{Enu}
\bfE(\nu) &=& \left( \bfepsilon_{\nu}, \bfepsilon_{s+\nu}, \ldots, \bfepsilon_{(N-1)s+\nu} \right),\\
\label{Xnu}
\bfX(\nu) &=& \left( \bfX_0(\nu), \ldots, \bfX_{N-1}(\nu) \right),
\end{eqnarray}
be
$d \times N$,
$d \times N$
and
$\{ d p(\nu) \} \times N$
random matrices,
where
\[
\bfX_n(\nu) = (\bfY_{ns+\nu-1}^{\top},\ldots,\bfY_{ns+\nu-p(\nu)}^{\top})^{\top},\]

$n=0,1,\ldots,N-1,$
denote
$\{ d p(\nu) \} \times 1$
random vectors.
The PVAR model can be reformulated as:
\begin{equation}
\label{modelZnu}
 \bfZ(\nu) = \bfB(\nu) \bfX(\nu) + \bfE(\nu), \; \nu = 1,\ldots,s,
\end{equation}
where the model parameters are collected in the $d \times \{ dp(\nu) \}$ matrix
$\bfB(\nu)$ which is defined as:
\begin{equation}
\label{Bnu}
\bfB(\nu) = \left( \bfPhi_1(\nu), \ldots, \bfPhi_{p(\nu)}(\nu) \right).
\end{equation}
Vectorizing, we obtain:
\begin{eqnarray}
 \bfz(\nu) &=& \{ \bfX^{\top}(\nu) \otimes \bfI_d \} \vec\{ \bfB(\nu) \} + \vec\{ \bfE(\nu) \}, \nonumber \\
           &=& \{ \bfX^{\top}(\nu) \otimes \bfI_d \} \bfbeta(\nu) + \bfe(\nu), \nonumber \\
\label{enu}
           &=& \{ \bfX^{\top}(\nu) \otimes \bfI_d \} \{ \bfR(\nu)\bfxi(\nu) + \bfb(\nu) \} + \bfe(\nu),
\end{eqnarray}
where
$\bfz(\nu) = \vec\{ \bfZ(\nu) \}$,
$\bfbeta(\nu) = \vec\{ \bfB(\nu) \}$,
$\bfe(\nu) = \vec\{ \bfE(\nu) \}$.

The covariance matrix of the random vector
$\bfe(\nu)$
is
$\bfI_N \otimes \bfSigma_{\bfepsilon}(\nu)$.

The multivariate least squares estimators of
$\bfxi(\nu)$, $\nu = 1,\ldots,s$ are obtained by minimizing
the generalized least squares criterion:
\begin{equation}
\label{GLS}
  S_G(\bfxi) = \sum_{\nu=1}^s \bfe^{\top}(\nu) \{ \bfI_N \otimes \bfSigma_{\bfepsilon}(\nu) \}^{-1}  \bfe(\nu),
\end{equation}
where
$\bfxi = (\bfxi^{\top}(1),\ldots,\bfxi^{\top}(s))^{\top}$
represents a
$\{ \sum_{\nu=1}^s K(\nu) \} \times 1$
vector. It may be worth nothing to mention that the GLS and LS estimation in a multiple equation model are identical if the regressors in all equations are the same (see for example a result for VAR models in~\citet[p.71]{Lu05}).
In the next subsections, we discuss separately the unrestricted and restricted cases.\\

\subsection{Unconstrained least squares estimators}\label{LSnoncont}
The least squares estimators are obtained by minimizing the ordinary least squares:
\begin{equation}
\label{ols}
  S(\bfbeta) = \sum_{\nu=1}^s \bfe^{\top}(\nu)   \bfe(\nu),
\end{equation}
where
$\bfbeta = (\bfbeta^{\top}(1),\ldots,\bfbeta^{\top}(s))^{\top}$
is the
$\{ d^2 \sum_{\nu=1}^s p(\nu) \} \times 1$
vector of model parameters.
To obtain the least squares estimators,
we differentiate
$S(\bfbeta)$
with respect to each parameter
$\bfPhi_k(\nu)$, $k=1,\ldots,p(\nu)$, $\nu=1,\ldots,s$.
Thus we obtain easily:
\[
  \frac{\partial S(\bfbeta)}{\partial \vec\{\bfPhi_k(\nu) \}} =
  -2\sum_{n=0}^{N-1} (\bfY_{ns+\nu-k} \otimes \bfepsilon_{ns+\nu}), \; k=1,\ldots,p(\nu), \; \nu=1,\ldots,s.
\]
Setting the derivatives equal to zero,
$k=1,\ldots,p(\nu)$,
gives the following system for a given season $\nu$:
\[
 \sum_{n=0}^{N-1} \left\{ \bfX_n(\nu) \otimes \bfepsilon_{ns+\nu} \right\} = \bfzero,
\]
where
$\bfzero$
is the
$\{ d^2 p(\nu) \} \times 1$
null vector.
Since
$\bfepsilon_{ns+\nu} = \bfY_{ns+\nu} - \{ \bfX_n^{\top}(\nu) \otimes \bfI_d \} \bfbeta(\nu)$,
the normal equations at season $\nu$ are:
\[
 \sum_{n=0}^{N-1} \{ \bfX_n(\nu) \otimes \bfY_{ns+\nu} \} =
 \left[ \sum_{n=0}^{N-1} \left\{ \bfX_n(\nu)\bfX_n^{\top}(\nu) \otimes \bfI_d \right\} \right] \bfbeta(\nu).
\]
Consequently,
the least squares estimators of
$\bfbeta(\nu)$
satisfy the relation:
\[
  \hat{\bfbeta}(\nu) = \left[ \{ \bfX(\nu)\bfX^{\top}(\nu) \}^{-1}\bfX(\nu) \otimes \bfI_d \right]\bfz(\nu),
\]
and the residuals are
$\hat{\bfepsilon}_{ns+\nu} = \bfY_{ns+\nu} - \{ \bfX_n^{\top}(\nu) \otimes \bfI_d \} \hat{\bfbeta}(\nu)$.
Using the
properties of the
$\vec(\cdot)$
operator,
it should be noted that an alternative expression for the least squares estimators is given by:
\begin{equation}
\label{hatBnu}
  \hat{\bfB}(\nu) = \bfZ(\nu) \bfX^{\top}(\nu) \{ \bfX(\nu)\bfX^{\top}(\nu) \}^{-1}.
\end{equation}
The asymptotic properties of the least squares estimators
in the unrestricted case are stated in Theorem~\ref{thm}.
The symbols
'$\cd$', '$\cp$'
and
'$\cas$'
stand for convergence in distribution, in probability and almost surely,
respectively,
and
$N_d(\bfmu, \bfSigma)$
denotes a $d$-dimensional normal distribution with mean $\bfmu$ and covariance matrix $\bfSigma$.

\begin{theorem}
\label{thm}
Let a time series be generated by equation~(\ref{pvar}).
Under the assumptions {\bf (A0)}, {\bf (A1)},  {\bf (A2)} and  {\bf (A3)}, for $\nu=1,\ldots,s$, we have
\begin{eqnarray}
\label{th1a}
  N^{-1/2} \sum_{n=0}^{N-1} \vec\{ \bfepsilon_{ns+\nu} \bfX_n^{\top}(\nu) \} &\cd&
  N_{d^2p(\nu)}\left( \bfzero,\bfPsi(\nu)
  \right), \\ \nonumber
  \bfPsi(\nu)&=& \sum_{h=-\infty}^{\infty}\mathbb{E}\left(\bfX_n(\nu)\bfX_{n-h}^{\top}(\nu) \otimes \bfepsilon_{ns+\nu}\bfepsilon_{(n-h)s+\nu}^{\top}\right)\\
\label{th1b}
  \hat{\bfbeta}(\nu) &\cas& \bfbeta(\nu), \\
\label{th1c}
  N^{1/2}\{ \hat{\bfbeta}(\nu) - \bfbeta(\nu) \}
  &\cd& N_{d^2p(\nu)}\left(\bfzero, \bfTheta (\nu) \right),
  \\ \nonumber
  \bfTheta(\nu)&=&\left(\bfOmega^{-1}(\nu) \otimes \bfI_d\right)\bfPsi(\nu)\left(\bfOmega^{-1}(\nu) \otimes \bfI_d\right)
\end{eqnarray}
where
$\bfOmega(\nu)$
corresponds to the
$\{ d p(\nu) \} \times \{ d p(\nu) \}$
covariance matrix of the
$\{ d p(\nu) \} \times 1$
random vector
$\bfX_n(\nu)$.
Furthermore, we also hawe
\begin{eqnarray}\nonumber
\label{th1d}
  N^{1/2}\{ \hat{\bfbeta} - \bfbeta\}
  &\cd& N_{sd^2p(\nu)}\left(\bfzero, \bfTheta \right),
\end{eqnarray}
where the asymptotic covariance matrix $\bfTheta$ is a block matrix, with the asymptotic variances given by $\bfTheta(\nu)$, $\nu=1,\dots,s$, and the asymptotic covariances given by:
\[\lim_{N\to\infty}\text{cov}\left(N^{1/2}\{ \hat{\bfbeta}(\nu) - \bfbeta(\nu) \},N^{1/2}\{ \hat{\bfbeta}(\nu') - \bfbeta(\nu') \}\right)=\left(\bfOmega^{-1}(\nu) \otimes \bfI_d\right)\sum_{h=-\infty}^{\infty}\mathbb{E}\left(\bfX_n(\nu)\bfX_{n-h}^{\top}(\nu') \otimes \bfepsilon_{ns+\nu}\bfepsilon_{(n-h)s+\nu'}^{\top}\right)\left(\bfOmega^{-1}(\nu') \otimes \bfI_d\right),\]
for $\nu \neq \nu'$ and $\nu, \nu' = 1,\ldots,s$.
%Furthermore,
%$N^{1/2}\{ \hat{\bfbeta}(\nu) - \bfbeta(\nu) \}$
%and
%$N^{1/2}\{ \hat{\bfbeta}(\nu') - \bfbeta(\nu') \}$
%are asymptotically independent,
%$\nu \neq \nu'$,
%$\nu, \nu' = 1,\ldots,s$.
\end{theorem}
The proof of Theorem \ref{thm} is postponed to Section \ref{app}.
\begin{remark} When $s=1$, we retrieve the well-known result on weak VAR obtained by \cite{FR07}.
\end{remark}
\begin{remark} If the moving average orders are null and when $d=1$ we retrieve the results obtained on weak periodic autoregressive model by~\citep{frs11}.
\end{remark}
\begin{remark}
\label{rem1}
In the standard strong PVAR case, i.e. when {\bf (A2)} is replaced by the assumption that   $\left(\boldsymbol{\epsilon}^{\ast}_n\right)_{n\in \mathbb{Z}}$ is  an iid sequence, we
have $$\bfPsi (\nu)=\bfOmega(\nu) \otimes \bfSigma_{\bfepsilon}(\nu).$$
Thus the asymptotic covariance matrix is reduced as
\[\bfTheta_S(\nu):=
  \{ \bfOmega^{-1}(\nu) \otimes \bfI_d \} \{ \bfOmega(\nu) \otimes \bfSigma_{\bfepsilon}(\nu) \}
  \{ \bfOmega^{-1}(\nu) \otimes \bfI_d \} =  \bfOmega^{-1}(\nu) \otimes \bfSigma_{\bfepsilon}(\nu) ,
\]
and we obtain the result of \cite{UD09}.

\noindent Generally, when the noise is not an independent sequence, this simplification can not be made and
we have $\bfPsi (\nu)\neq\bfOmega(\nu) \otimes \bfSigma_{\bfepsilon}(\nu)$.
The true asymptotic covariance matrix $\bfTheta(\nu)$ obtained
in the weak PVAR framework can be very different from $\bfTheta_S(\nu)$.
As a consequence, for the statistical inference on the
parameter,  the ready-made software used to
fit PVAR do not provide a correct estimation of $\bfTheta(\nu)$ for weak PVAR processes because the standard
time series analysis software use empirical estimators of $\bfTheta_S(\nu)$.
The problem also holds in the weak PARMA case (see
\cite{frs11} and the references therein). This is why it is interesting to find  an estimator of $\bfTheta(\nu)$  which is consistent
for both weak and strong PVAR cases.
\end{remark}

%--------------------------------------------------------------------------------------------------------------------------------------
\subsection{Least squares estimation with linear constraints on the parameters}\label{LScont}
\noindent When the parameters satisfy the linear constraint~(\ref{constraints}), the least squares estimators of $\bfxi(\nu)$, $\nu = 1,\ldots,s$, minimize the generalized criterion~(\ref{GLS}), which is not equivalent to~(\ref{ols}), see~\citet{Lu05}, amongst others. Recall that from~(\ref{enu}) we have the following relation:
\[
  \bfe(\nu) = \bfz(\nu) - \{ \bfX^{\top}(\nu) \otimes \bfI_d \} \{ \bfR(\nu) \bfxi(\nu) + \bfb(\nu) \},
\]
which is convenient to derive the asymptotic properties of the least squares estimator of $\bfxi(\nu)$.

\noindent Proceeding as in the previous section, it is possible to show that the least squares estimator $\hat{\bfxi}(\nu)$ of $\bfxi(\nu)$ is given by:
\begin{align*}
\hat{\bfxi}(\nu) &=
  \left[ \bfR^{\top}(\nu) \{ \bfX(\nu)\bfX^{\top}(\nu) \otimes \bfSigma^{-1}_{\bfepsilon}(\nu) \} \bfR(\nu) \right]^{-1}
  \bfR^{\top}(\nu) \{ \bfX(\nu)  \otimes \bfSigma^{-1}_{\bfepsilon}(\nu) \}
 \\& \hspace{7cm}\times  \left[\bfz(\nu) - \{ \bfX^{\top}(\nu) \otimes \bfI_d \}\bfb(\nu) \right].
\end{align*}
Furthermore, the following relation is satisfied:
\begin{align*}
  N^{1/2} &\{  \hat{\bfxi}(\nu) - \bfxi(\nu) \} \\&=
  N^{1/2} \left[ \bfR^{\top}(\nu) \{ \bfX(\nu)\bfX^{\top}(\nu) \otimes \bfSigma^{-1}_{\bfepsilon}(\nu) \} \bfR(\nu) \right]^{-1}
  \bfR^{\top}(\nu) \{ \bfX(\nu)  \otimes \bfSigma^{-1}_{\bfepsilon}(\nu) \} \bfe(\nu)\\
  &= \left[ \bfR^{\top}(\nu)\frac{1}{N} \{ \bfX(\nu)\bfX^{\top}(\nu) \otimes \bfSigma^{-1}_{\bfepsilon}(\nu) \} \bfR(\nu) \right]^{-1}
  \bfR^{\top}(\nu)\\& \hspace{6cm}\times \{ \bfI_{dp(\nu)}  \otimes \bfSigma^{-1}_{\bfepsilon}(\nu) \} N^{-1/2} \vec\{ \bfE(\nu) \bfX^\top(\nu)  \}.
\end{align*}
Consequently, under the conditions of Theorem~\ref{thm}, the estimator $\hat{\bfxi}(\nu)$
is consistent for $\bfxi(\nu)$, and $\hat{\bfxi}(\nu)$ follows asymptotically a normal distribution, that is:
\begin{align}
\label{NAcont}
N^{1/2} \{  \hat{\bfxi}(\nu) - \bfxi(\nu) \} \cd
  N_{K(\nu)}\left(\bfzero,
              \bfTheta^{\bfxi}(\nu)
            \right),
\end{align}
where
\begin{align*}
  \bfTheta^{\bfxi}(\nu)&=\left[ \bfR^{\top}(\nu) \{ \bfOmega(\nu) \otimes \bfSigma^{-1}_{\bfepsilon}(\nu) \} \bfR(\nu) \right]^{-1}
\bfR^{\top}(\nu) \{ \bfI_{dp(\nu)} \otimes \bfSigma^{-1}_{\bfepsilon}(\nu) \}  \bfPsi(\nu) \\&\qquad\times\{ \bfI_{dp(\nu)} \otimes \bfSigma^{-1}_{\bfepsilon}(\nu) \}   \bfR(\nu)
 \left( \left[ \bfR^{\top}(\nu) \{ \bfOmega(\nu) \otimes \bfSigma^{-1}_{\bfepsilon}(\nu) \} \bfR(\nu) \right]^{-1}\right)^{\top}.
\end{align*}
Moreover we have
\begin{align*}
&\{ \bfI_{dp(\nu)} \otimes \bfSigma^{-1}_{\bfepsilon}(\nu) \}  \bfPsi(\nu) \{ \bfI_{dp(\nu)} \otimes \bfSigma^{-1}_{\bfepsilon}(\nu) \}
  \\&=\left(\bfI_{dp(\nu)} \otimes \bfSigma^{-1}_{\bfepsilon}(\nu)\right)\sum_{h=-\infty}^{\infty}\mathbb{E}\left(\bfX_n(\nu)\bfX_{n-h}^{\top}(\nu) \otimes \bfepsilon_{ns+\nu}\bfepsilon_{(n-h)s+\nu}^{\top}\right)\left(\bfI_{dp(\nu)} \otimes \bfSigma^{-1}_{\bfepsilon}(\nu)\right)
\\
 &= \sum_{h=-\infty}^{\infty}\mathbb{E}\left[\bfX_n(\nu)\bfX_{n-h}^{\top}(\nu) \otimes \bfSigma^{-1}_{\bfepsilon}(\nu)\bfepsilon_{ns+\nu}\bfepsilon_{(n-h)s+\nu}^{\top}\bfSigma^{-1}_{\bfepsilon}(\nu)\right].
\end{align*}
It should be noted that the estimator $\hat{\bfxi}(\nu)$ is unfeasible in practice, since it relies on the unknown matrix $\bfSigma_{\bfepsilon}(\nu)$. A feasible estimator is given by:
\begin{align*}
\hat{\hat{\bfxi}}(\nu) &=
  \left[ \bfR^{\top}(\nu) \{ \bfX(\nu)\bfX^{\top}(\nu) \otimes \tilde{\bfSigma}^{-1}_{\bfepsilon}(\nu) \} \bfR(\nu) \right]^{-1}
  \bfR(\nu) \{ \bfX(\nu)  \otimes \tilde{\bfSigma}^{-1}_{\bfepsilon}(\nu) \}
\\& \hspace{7cm}\times   [\bfz(\nu) - \{ \bfX^{\top}(\nu) \otimes \bfI_d \}\bfb(\nu)],
\end{align*}
\noindent where $\tilde{\bfSigma}_{\bfepsilon}(\nu)$ denotes a consistent estimator of the covariance matrix $\bfSigma_{\bfepsilon}(\nu)$ for $\nu=1,\ldots,s$.
A possible candidate is obtained from the unconstrained least squares estimators:
\[
\tilde{\bfSigma}_{\bfepsilon}(\nu) = \{ N-d p(\nu) \}^{-1} \left\{ \bfZ(\nu) - \hat{\bfB}(\nu)\bfX(\nu) \right\}
                                     \left\{ \bfZ(\nu) - \hat{\bfB}(\nu)\bfX(\nu) \right\}^{\top},
\]
where $\hat{\bfB}(\nu)$ represents the unconstrained least squares estimators~(\ref{hatBnu}) obtained in Section~\ref{LSnoncont}. The resulting estimator of
$\bfbeta(\nu)$ is given by
$\hat{\hat{\bfbeta}}(\nu) = \bfR(\nu)\hat{\hat{\bfxi}}(\nu) + \bfb(\nu)$, and its asymptotic distribution  is normal:
\[
  N^{1/2} \{  \hat{\hat{\bfbeta}}(\nu) - \bfbeta(\nu) \} \cd
  N_{d^2p(\nu)}\left(\bfzero,
             \bfR(\nu)
              \bfTheta^{\bfxi}(\nu)
             \bfR^{\top}(\nu)
            \right).
\]
The proof of the above result follows, using arguments similar to those of Theorem~\ref{thm}.
\begin{remark}\label{rem2}
In the standard strong PVAR case, i.e. when {\bf (A2)} is replaced by the assumption that   $\left(\boldsymbol{\epsilon}^{\ast}_n\right)_{n\in \mathbb{Z}}$ is  an iid sequence and in view of Remark~\ref{rem1}, we
have
\begin{align*}
  N^{1/2} \{  \hat{\bfxi}(\nu) - \bfxi(\nu) \}& \cd
  N_{K(\nu)}\left(\bfzero,\bfTheta_S^{\bfxi}(\nu)=:
             \left[ \bfR^{\top}(\nu) \{ \bfOmega(\nu) \otimes \bfSigma^{-1}_{\bfepsilon}(\nu) \} \bfR(\nu) \right]^{-1}
            \right),
\\
  N^{1/2} \{  \hat{\hat{\bfbeta}}(\nu) - \bfbeta(\nu) \} &\cd
  N_{d^2p(\nu)}\left(\bfzero,
             \bfR(\nu)
             \bfTheta_S^{\bfxi}(\nu)
             \bfR^{\top}(\nu)
            \right),
\end{align*}
which are the results obtained by \cite{UD09}.
\end{remark}
\subsection{Example of analytic computation of $\bfTheta(\nu)$ and $\bfTheta_S(\nu)$}\label{exple}
\noindent Consider a bi-variate periodic white noise  defined by:
\begin{eqnarray}\label{weaknoise}
\bfepsilon_{ns+\nu} =
\begin{pmatrix}
\bfepsilon_{1,ns+\nu} \\
\bfepsilon_{2,ns+\nu}
\end{pmatrix}=\bfM_{\nu}^\top
\begin{pmatrix}
\eta_{1,ns+\nu}\eta_{1,ns+\nu-1}\eta_{1,ns+\nu-2}\cdots\eta_{1,ns+\nu-m} \\
\eta_{2,ns+\nu}\eta_{2,ns+\nu-1}\eta_{2,ns+\nu-2}\cdots\eta_{2,ns+\nu-m}
\end{pmatrix}, %\quad \nu=1,2;
\end{eqnarray}
where $m>0$ is a fixed integer,  $\bfeta_{t} =(\eta_{1,t},\eta_{2,t})^\top$ iid $\mathcal{N}(\bfzero,\bfI_2)$ and $\bfM_{\nu}$ is the upper triangular matrix satisfying the equation $\bfM^{T}_{\nu} \bfM_{\nu} = \bfSigma_{\bfepsilon}(\nu)$. The periodic process $\bfepsilon_{ns+\nu}$ is a weak white noise because $\mathbb{E}[\bfepsilon_t] = \bfzero$ for all $t$, $\mathbb{E}[\bfepsilon_t\bfepsilon_{t'}^\top] =\bfzero$ for all $t\neq t'$, $\mathbb{E}[\bfepsilon_{ns+\nu}\bfepsilon_{ns+\nu}^\top] = \bfSigma_{\bfepsilon}(\nu)$.
The variables $\bfepsilon_t$ and $\bfepsilon_{t'}$ are dependent if $|t-t'|\leq m$ but they are independent for $|t-t'| > m$. The process \eqref{weaknoise} can be viewed as a multivariate extension of univariate weak noises considered in~\citet{frs11}.

\noindent The results of Section \ref{result1}  is particularized in the following PVAR case of order one with $s=2$ of the form:
\begin{equation}
\label{PAR1}
\left\{\begin{array}{ll}
\bfY_{2n+1}=\bfPhi(1)\bfY_{2n}+\bfepsilon_{2n+1}\\ %\vspace{0.1cm}\\
\bfY_{2n+2}=\bfPhi(2)\bfY_{2n+1}+\bfepsilon_{2n+2}
\end{array}\right.,
\end{equation}
where the  unknown parameter is $\bfbeta(\nu)=\vec{(\bfPhi(\nu))}$ and the innovation process $(\bfepsilon_{2n+\nu})_{n\in\mathbb{Z}}$  is  given by \eqref{weaknoise}. For simplification, we assume that in \eqref{PAR1}, $\bfPhi(\nu)$ and $\bfSigma_{\nu}$ are diagonals:
\begin{eqnarray*}
\bfPhi(\nu)=\begin{pmatrix}
\phi_{11}(\nu)&0\\
0&\phi_{22}(\nu)
\end{pmatrix}
\qquad\text{and}\qquad\bfSigma_{\bfepsilon}(\nu)=
\begin{pmatrix}
\sigma_{11}(\nu)&0\\
0&\sigma_{22}(\nu)
\end{pmatrix}.
\end{eqnarray*}
From \eqref{VARrepre} and under {\bf (A1)} we deduce that:
\begin{equation}
\label{PAR1linear}
\left\{\begin{array}{ll}
\bfY_{2n+1}=\bfepsilon_{2n+1}+{\bfPhi^{-1}(2)}\sum_{i\geq1}\bfPhi(2)^i\bfPhi^i(1)\bfepsilon_{2(n-i)+2}
\\ %\vspace{0.1cm}\\
\bfY_{2n+2}=\bfPhi(2)\bfepsilon_{2n+1}+\sum_{i\geq0}\bfPhi^i(2)\bfPhi^i(1)\bfepsilon_{2(n-i)+2}
\end{array}\right..
\end{equation}
With our notations $\bfX_n(\nu) = \bfY_{2n+\nu-1}$ and by using \eqref{PAR1linear}, it follows that:
\begin{align}\label{omegaex}
\bfOmega(1)&=\mathbb{E}\left(\bfY_{2n}\bfY_{2n}^\top\right)=
\begin{pmatrix}\bfOmega_{11}(1)&0\\0 &\bfOmega_{22}(1)
 \end{pmatrix}\quad \text{ and }\quad \bfOmega(2)=\mathbb{E}\left(\bfY_{2n+1}\bfY_{2n+1}^\top\right)=
\begin{pmatrix}\bfOmega_{11}(2)&0\\0 &\bfOmega_{22}(2)
 \end{pmatrix},
\end{align}
where for $i=1,2$ we have
%\begin{align*}
%\bfOmega_{11}(1)&=\frac{\phi_{11}^2(2)\sigma_{11}(1)\left(1-\phi_{11}^2(1)\phi_{11}^2(2)\right)+\sigma_{11}(2)}{1-\phi_{11}^2(1)\phi_{11}^2(2)},\qquad
%\bfOmega_{22}(1)=\frac{\phi_{22}^2(2)\sigma_{22}(1)\left(1-\phi_{22}^2(1)\phi_{22}^2(2)\right)+\sigma_{22}(2)}{1-\phi_{22}^2(1)\phi_{22}^2(2)},\\
%\bfOmega_{11}(2)&=\frac{\phi_{11}^2(2)\sigma_{11}(1)\left(1-\phi_{11}^2(1)\phi_{11}^2(2)\right)+\sigma_{11}(2)\phi_{11}^2(1)\phi_{11}^2(2)}{\phi_{11}^2(2)\left( 1-\phi_{11}^2(1)\phi_{11}^2(2)\right)} \text{ and }
%\bfOmega_{22}(2)=\frac{\phi_{22}^2(2)\sigma_{22}(1)\left(1-\phi_{22}^2(1)\phi_{22}^2(2)\right)+\sigma_{22}(2)\phi_{22}^2(1)\phi_{22}^2(2)}{\phi_{22}^2(2)\left( 1-\phi_{22}^2(1)\phi_{22}^2(2)\right)},
%\end{align*}
\begin{align}\label{omegaii}
\bfOmega_{ii}(1)&=\frac{\phi_{ii}^2(2)\sigma_{ii}(1)\left(1-\phi_{ii}^2(1)\phi_{ii}^2(2)\right)+\sigma_{ii}(2)}{1-\phi_{ii}^2(1)\phi_{ii}^2(2)}, \text{ and }
\bfOmega_{ii}(2)=\frac{\phi_{ii}^2(2)\sigma_{ii}(1)\left(1-\phi_{ii}^2(1)\phi_{ii}^2(2)\right)+\sigma_{ii}(2)\phi_{ii}^2(1)\phi_{ii}^2(2)}{\phi_{ii}^2(2)\left( 1-\phi_{ii}^2(1)\phi_{ii}^2(2)\right)}.
\end{align}
In view of Remark \ref{rem1} and using \eqref{omegaex}, a simple calculation implies that
\begin{equation}\label{invJ_2}
\bfTheta_S(\nu) =  \bfOmega^{-1}(\nu) \otimes \bfSigma_{\bfepsilon}(\nu)=\begin{pmatrix} \bfOmega_{11}^{-1}(\nu) \bfSigma_{\bfepsilon}(\nu)  &  0 \vspace{0.2cm}\\  0 &  \bfOmega_{22}^{-1}(\nu)  \bfSigma_{\bfepsilon}(\nu)\end{pmatrix},\quad\nu=1,2,
\end{equation}
where $\bfOmega_{ii}(1)$, $\bfOmega_{ii}(2)$ are given by \eqref{omegaii} and we obtain the result of \cite{UD09}.

\noindent We now investigate a similar tractable expression for $ \bfPsi(\nu)= \sum_{h=-\infty}^{\infty}\mathbb{E}\left(\bfX_n(\nu)\bfX_{n-h}^{\top}(\nu) \otimes \bfepsilon_{ns+\nu}\bfepsilon_{(n-h)s+\nu}^{\top}\right)$. Using \eqref{weaknoise} and \eqref{PAR1linear}, the matrix $\bfPsi(\nu)$ is given  by
\begin{align}\label{I_1}
\bfPsi(1) =  \mbox{Diag}\begin{pmatrix}\bfPsi_{11}(1), \bfPsi_{22}(1),\bfPsi_{33}(1),\bfPsi_{44}(1)
\end{pmatrix}\quad\text{and}\quad \bfPsi(2) =  \mbox{Diag}\begin{pmatrix}\bfPsi_{11}(2), \bfPsi_{22}(2),\bfPsi_{33}(2),\bfPsi_{44}(2)
\end{pmatrix},
\end{align}
where
\begin{align*}
\bfPsi_{11}(1) &=3^{m-1}\phi_{11}^2(2)\sigma_{11}^2(1)+\sigma_{11}(1)\sigma_{11}(2)\left(\sum_{i=0}^{\lfloor \frac{m-1}{2} \rfloor}3^{m-2i}\phi_{11}^{2i}(1)\phi_{11}^{2i}(2)+
\frac{\phi_{11}^{2\left(\lfloor \frac{m-1}{2} \rfloor+1\right)}(1)\phi_{11}^{2\left(\lfloor \frac{m-1}{2} \rfloor+1\right)}(2)}{1-\phi_{11}^{2}(1)\phi_{11}^{2}(2)}\right),\\
\bfPsi_{22}(1) &=\phi_{11}^2(2)\sigma_{11}(1)\sigma_{22}(1)+
\frac{\sigma_{22}(1)\sigma_{11}(2)}{1-\phi_{11}^{2}(1)\phi_{11}^{2}(2)},\quad
\bfPsi_{33}(1) =\phi_{22}^2(2)\sigma_{11}(1)\sigma_{22}(1) +\frac{\sigma_{11}(1)\sigma_{22}(2)}{1-\phi_{22}^{2}(1)\phi_{22}^{2}(2)},\\
\bfPsi_{44}(1) &=3^{m-1}\phi_{22}^2(2)\sigma_{22}^2(1)+
\sigma_{22}(1)\sigma_{22}(2)\left(\sum_{i=0}^{\lfloor \frac{m-1}{2} \rfloor}3^{m-2i}\phi_{22}^{2i}(1)\phi_{22}^{2i}(2)+
\frac{\phi_{22}^{2\left(\lfloor \frac{m-1}{2} \rfloor+1\right)}(1)\phi_{22}^{2\left(\lfloor \frac{m-1}{2} \rfloor+1\right)}(2)}{1-\phi_{22}^{2}(1)\phi_{22}^{2}(2)}\right),\\
\bfPsi_{11}(2) &=3^{m}\sigma_{11}^2(1)+\frac{\sigma_{11}^{2}(2)}{\phi_{11}^{2}(2)}\left(\sum_{i=1}^{\lfloor \frac{m}{2} \rfloor}3^{m-2i+1}\phi_{11}^{2i}(1)\phi_{11}^{2i}(2)+
\frac{\phi_{11}^{2\left(\lfloor \frac{m}{2} \rfloor+1\right)}(1)\phi_{11}^{2\left(\lfloor \frac{m}{2} \rfloor+1\right)}(2)}{1-\phi_{11}^{2}(1)\phi_{11}^{2}(2)}\right),\\
\bfPsi_{22}(2) &=\sigma_{11}(1)\sigma_{22}(1)+\sigma_{22}(2)\sigma_{11}(2)
\frac{\phi_{11}^{2}(1)}{1-\phi_{11}^{2}(1)\phi_{11}^{2}(2)},\quad
\bfPsi_{33}(2) =\sigma_{11}(1)\sigma_{22}(1) +\sigma_{11}(2)\sigma_{22}(2)\frac{\phi_{22}^{2}(1)}{1-\phi_{22}^{2}(1)\phi_{22}^{2}(2)},\\
\bfPsi_{44}(2) &=3^{m}\sigma_{22}^2(1)+\frac{\sigma_{22}^{2}(2)}{\phi_{22}^{2}(2)}
\left(\sum_{i=1}^{\lfloor \frac{m}{2} \rfloor}3^{m-2i+1}\phi_{22}^{2i}(1)\phi_{22}^{2i}(2)+
\frac{\phi_{22}^{2\left(\lfloor \frac{m}{2} \rfloor+1\right)}(1)\phi_{22}^{2\left(\lfloor \frac{m}{2} \rfloor+1\right)}(2)}{1-\phi_{22}^{2}(1)\phi_{22}^{2}(2)}\right).
\end{align*}
\noindent From \eqref{I_1} we deduce that
\begin{eqnarray}
\label{sandwich}
 \bfTheta(\nu)=\left(\bfOmega^{-1}(\nu) \otimes \bfI_d\right)\bfPsi(\nu)\left(\bfOmega^{-1}(\nu) \otimes \bfI_d\right),\quad \nu=1,2,
\end{eqnarray}
where  $\bfOmega_{ii}(1)$, $\bfOmega_{ii}(2)$ are given by \eqref{omegaii} for $i=1,2$.

\noindent For instance when
\begin{eqnarray*}
\bfPhi(1)=\begin{pmatrix}
0.3&0.0\\
0.0&-0.6
\end{pmatrix},
\quad\bfPhi(2)=\begin{pmatrix}
-0.7&0.0\\
0.0&0.15
\end{pmatrix},
\quad
\bfSigma_{\bfepsilon}(1)=
\begin{pmatrix}
1.5&0.0\\
0.0&2.5
\end{pmatrix}
\quad\text{and}\quad\bfSigma_{\bfepsilon}(2)=
\begin{pmatrix}
1&0.0\\
0.0&0.5
\end{pmatrix}
\end{eqnarray*}
we have
\begin{eqnarray*}
\bfTheta_S(1)=\begin{pmatrix}
0.84 & 0.00& 0.00& 0.00\\
0.00 & 1.40& 0.00& 0.00\\
 0.00 & 0.00 &2.68& 0.00\\
0.00&  0.00 &0.00& 4.46
\end{pmatrix},
\quad\bfTheta_S(2)=\begin{pmatrix}
0.69 &0.00 &0.00 &0.00\\
0.00 &0.34 &0.00& 0.00\\
0.00& 0.00 &0.38& 0.00\\
0.00& 0.00& 0.00& 0.19
\end{pmatrix},
\end{eqnarray*}
and the analytic computation of $\bfTheta(\nu)$ is given for  $m=1,2$ in the following Table. %\ref{tab1}.
%\begin{table}[H]
\begin{center}
%\scriptsize
\begin{tabular}{|c|c|c|}
\hline \hline
$\nu$&  \multicolumn{2}{c}{$\bfTheta(\nu)$} \\
\hline
 &$m=1$ &$m=2$ \\
\cline{1-3}
$\nu=1$& $\begin{pmatrix} 1.79&  0.00 &0.00 & 0.00\\ 0.00  &1.40& 0.00 & 0.00\\ 0.00 & 0.00& 2.68&  0.00\\ 0.00&  0.00& 0.00& 12.42
\end{pmatrix}$ & $\begin{pmatrix} 2.48 & 0.00& 0.00 & 0.00\\0.00&  1.40& 0.00 & 0.00\\0.00&  0.00& 2.68 & 0.00\\0.00 & 0.00& 0.00&13.32\end{pmatrix}$ \\%\rule[7pt]{0pt}{20pt}\\
\hline
$\nu=2$& $\begin{pmatrix} 3.23& 0.00 &0.00& 0.00\\ 0.00& 1.79& 0.00& 0.00\\ 0.00& 0.00& 0.56& 0.00\\ 0.00& 0.00& 0.00& 2.71
\end{pmatrix}$ & $\begin{pmatrix} 9.72& 0.00& 0.00& 0.00\\ 0.00& 1.79& 0.00& 0.00\\ 0.00 &0.00 &0.56& 0.00\\ 0.00 &0.00 &0.00& 8.13
\end{pmatrix}$ \\%\rule[7pt]{0pt}{20pt}\\
\hline
\end{tabular}
\end{center}
%\label{tab1}
%\end{table}
%\normalsize
\noindent
For this particular weak noise, we draw the conclusion that the discrepancy between the two matrices $\bfTheta(\nu)$ and $\bfTheta_S(\nu)$ is important even for small $m$. For statistical inference problem, including in particular, the significance tests on the parameters, the assumption of independent errors can be quite misleading when analyzing data from PVAR models with dependent errors.
Thus the standard methodology needs however to be adapted to take into account the possible lack of independence of the error terms.
\section{Estimating the asymptotic variance matrix}\label{estimOmega}\label{result2}
\noindent For statistical inference problem, the asymptotic variance $\bfTheta(\nu)$ has to be estimated. In particular
Theorem \ref{thm} can be used to obtain confidence intervals and significance tests for the parameters.
\subsection{Estimation of the asymptotic matrix $\bfOmega(\nu)$}\label{i}
\noindent The matrix $\bfOmega(\nu)$ can be  estimated empirically by the square matrix $\hat{\bfOmega}_N(\nu)$ of order $dp(\nu)$ defined by:
\begin{equation}\label{estim_J}
\hat{\bfOmega}_N(\nu)=\frac{1}{N} \bfX_n(\nu) \bfX_n^{\top}(\nu)  .
\end{equation}
The convergence of $\hat{\bfOmega}_N(\nu)$ to $\bfOmega(\nu)$ is  proved in \eqref{omega}.

\noindent In the standard strong PVAR case, in view of remark~\ref{rem1}, we have  $\hat{\bfTheta}_S(\nu):= \hat{\bfOmega}_N^{-1}(\nu) \otimes \tilde{\bfSigma}_{\bfepsilon}(\nu)$. Thus $\hat{\bfTheta}_S(\nu)$  is a consistent estimator of $\bfTheta_S(\nu)$.
In the general weak PVAR case, this estimator is not consistent when $\bfPsi (\nu)\neq\bfOmega(\nu) \otimes \bfSigma_{\bfepsilon}(\nu)$. So we need a consistent estimator of $\bfPsi (\nu)$.
\subsection{Estimation of the asymptotic matrix $\bfPsi (\nu)$}\label{ii}
\noindent For all $n\in\mathbb{Z}$, let
\begin{align}\label{Ht_process}
\bfW_n(\nu):=\vec\{ \bfepsilon_{ns+\nu}\bfX_n^{\top}(\nu) \}.
\end{align}
We shall see in the proof of Theorem~\ref{thm} that
\begin{align}\label{matI}
\bfPsi(\nu)=\lim_{N\rightarrow\infty}\var\left( \frac{1}{\sqrt{N}}\sum_{n=0}^{N-1}\bfW_{n}(\nu)\right)=\sum_{h=-\infty}^{\infty}
\cov\left(\bfW_{n}(\nu),\bfW_{n-h}(\nu)\right) .
\end{align}
The estimation of the long-run variance (LRV) matrix $\bfPsi(\nu)$ is more complicated. In the literature, two types of estimators are generally
employed: heteroskedasticity and autocorrelation consistent (HAC) estimators based on kernel methods (see \cite{newey} and \cite{A91} for general references, and \cite{FZ07} for an application to testing strong linearity in weak ARMA models) and the spectral density estimators (see e.g. \citep{B74} and \cite{haan} for a general reference; see also \cite{BMF11} for an application to a weak VARMA model).
\subsubsection{Spectral density estimation of LRV matrix $\bfPsi(\nu)$}
\noindent Following the arguments developed in \cite{BMCF12}, the matrix $\bfPsi(\nu)$ can be estimated using Berk's approach (see \cite{B74}). More precisely, by interpreting $\bfPsi(\nu)/2\pi$ as the spectral density of the stationary process $(\bfW_n(\nu))_{n\in\mathbb{Z}}$ evaluated at frequency $0$, we can use a parametric autoregressive estimate of the spectral density of $(\bfW_n(\nu))_{n\in\mathbb{Z}}$ in order to estimate the matrix $\bfPsi(\nu)$.

\noindent The process $(\bfW_n(\nu))_{n\in\mathbb{Z}}$ is a measurable function of $\left\lbrace \bfepsilon_{ns+\nu-k},k\geq 0\right\rbrace $.
The stationary process $(\bfW_n(\nu))_{n\in\mathbb{Z}}$ admits the following Wold decomposition $\bfW_n(\nu)=u_{ns+\nu}+\sum_{k=1}^{\infty}\psi_k(\nu)u_{ns+\nu-k}$, where $(u_{ns+\nu})_{n\in\mathbb{Z}}$ is a $(d^2p(\nu))-$variate periodic weak white noise with variance matrix $\Sigma_u(\nu)$.

\noindent Assume that $\Sigma_u(\nu)$ is non-singular, that $\sum_{k=1}^{\infty}\left\|\psi_k(\nu)\right\|<\infty$, and that $\det(\bfI_{d^2p(\nu)}+\sum_{k=1}^{\infty}\psi_k(\nu)z^k)\neq 0$ if $\left|z\right|\leq 1$. Then $(\bfW_n(\nu))_{n\in\mathbb{Z}}$ admits a weak multivariate $\mathrm{AR}(\infty)$ representation (see \cite{A57}) of the form
\begin{equation}\label{AR_infty}
\Phi(L,\nu)\bfW_n(\nu):=\bfW_n(\nu)-\sum_{k=1}^{\infty}\Phi_k(\nu)\bfW_{n-k}(\nu)=u_{ns+\nu},
\end{equation}
such that $\sum_{k=1}^{\infty}\left\|\Phi_k(\nu)\right\|<\infty$ and $\det\left\lbrace \Phi(z,\nu)\right\rbrace \neq 0$ if $\left|z\right|\leq 1$. %

\noindent Thanks to the previous remarks, the estimation of $\bfPsi(\nu)$ is therefore based on the following expression
$$\bfPsi(\nu)=\Phi^{-1}(1,\nu)\Sigma_u(\nu)\Phi^{-1}(1,\nu).$$
Consider the regression of $\bfW_n(\nu)$ on $\bfW_{n-1}(\nu),\dots,\bfW_{n-r}(\nu)$ defined by
\begin{equation}\label{AR_tronquee}
\bfW_n(\nu)=\sum_{k=1}^{r}\Phi_{r,k}(\nu)\bfW_{n-k}(\nu)+u_{r,ns+\nu},
\end{equation}
where $u_{r,ns+\nu}$ is uncorrelated with $\bfW_{n-1}(\nu),\dots,\bfW_{n-r}(\nu)$.
Since $\bfW_n(\nu)$ is not observable, we introduce $\hat{\bfW}_n(\nu)\in\mathbb{R}^{d^2p(\nu)}$ obtained by replacing $\bfepsilon_{ns+\nu}$ by $\hat{\bfepsilon}_{ns+\nu}$ and $\beta(\nu)$ by $\hat{\beta}(\nu)$ in $\eqref{Ht_process}$:
\begin{align}\label{Ht_chapeau}
\hat{\bfW}_n(\nu)=\vec\{ \hat{\bfepsilon}_{ns+\nu}\bfX_n^{\top}(\nu) \} \ ,
\end{align}
where $\hat{\bfepsilon}_{ns+\nu}$ represents the unconstrained least squares residual.

\noindent Let $\hat{\Phi}_r(z,\nu)=\bfI_{dp(\nu)}-\sum_{k=1}^r\hat\Phi_{r,k}(\nu)z^k$, where $\hat\Phi_{r,1}(\nu),\dots,\hat\Phi_{r,r}(\nu)$ denote the coefficients of the LS regression
of $\hat{\bfW}_n(\nu)$ on $\hat{\bfW}_{n-1}(\nu),\dots,\hat{\bfW}_{n-r}(\nu)$. Let $\hat{u}_{r,ns+\nu}$ be the residuals of this regression and  let $\hat{\Sigma}_{\hat{u}_r}(\nu)$ be the empirical variance of $\hat{u}_{r,\nu},\dots,\hat{u}_{r,(N-1)s+\nu}$.

\noindent In the case of linear processes with independent innovations, \cite{B74} has shown that the spectral density can be consistently estimated by fitting autoregressive models of order $r=r(N)$, whenever $r$ tends to infinity and $r^3/N$ tends to $0$ as $N$ tends to infinity.  There are differences with \cite{B74}: $(\bfW_n(\nu))_{n\in\mathbb{Z}}$ is multivariate, is not directly observed and is replaced by $(\hat{\bfW}_n(\nu))_{n\in\mathbb{Z}}$. It is shown that this result remains valid for the multivariate linear process $(\bfW_n(\nu))_{n\in\mathbb{Z}}$ with non-independent innovations (see \cite{BMCF12,BMF11}, for references in weak (multivariate) ARMA models).

\noindent The asymptotic study of the estimator of $\bfPsi(\nu)$ using the spectral density method is given in the following theorem.
\begin{theorem}\label{convergence_Isp}
In addition to the assumptions of Theorem~\ref{thm}, assume that the process $(\bfW_n(\nu))_{t\in\mathbb{Z}}$ defined in \eqref{Ht_process} admits a periodic VAR$(\infty)$ representation \eqref{AR_infty}, where $\|\Phi_k(\nu)\|=\mathrm{o}(k^{-2})$ as $k\to\infty$, the roots of $\det(\Phi(z,\nu))=0$ are outside the unit disk, and $\Sigma_u(\nu)=\var(u_{ns+\nu})$ is non-singular. Moreover we assume that $\mathbb{E}\|\boldsymbol{\epsilon}^{\ast}_n\|^{8+4\kappa}<\infty$
%and $\sum_{k=0}^{\infty}\left\{\alpha_{\boldsymbol{\epsilon}}(k)\right\}^{\frac{\kappa}{2+\kappa}}<\infty$
for some $\kappa>0$.
Then, the spectral estimator of $\bfPsi(\nu)$:
$$\hat{\bfPsi}^{\mathrm{SP}}(\nu):=\hat{\Phi}_r^{-1}(1,\nu)\hat{\Sigma}_{\hat{u}_r}(\nu)\hat{\Phi}_r^{'-1}(1,\nu)\cp \bfPsi(\nu)=\Phi^{-1}(1,\nu)\Sigma_u(\nu)\Phi^{-1}(1,\nu)$$
when $r=r(N)\to\infty$ and $r^{3}/N\to0$ as $N\to\infty$.
\end{theorem}
\noindent The proof of this theorem is similar to that given by \cite[Theorem 3]{BMIA22} and it is omitted.
\subsubsection{HAC estimation of LRV matrix $\bfPsi(\nu)$}
\noindent Let
\begin{equation}\label{mom1}
\Lambda_h(\nu)=\cov\left(\bfW_{n}(\nu),\bfW_{n-h}(\nu)\right)=\mathbb{E}\left(\bfW_{n}(\nu)\bfW_{n-h}^\top(\nu)\right).
\end{equation}
The sum $\sum_{h=-\infty}^{\infty}\Lambda_h(\nu)$ is well defined (see the proof of Theorem~\ref{thm}).
From the stationarity of the centered process $(\bfW_n(\nu))_{n\in\mathbb{Z}}$ and by the Lebesgue theorem, we have
\begin{align}\label{I-noyau}
\bfPsi(\nu)&=\lim_{N\rightarrow\infty}\var\left\lbrace \frac{1}{\sqrt{N}}\sum_{n=0}^{N-1} \bfW_n(\nu) \right\rbrace =\lim_{N\rightarrow\infty}\frac{1}{N}\sum_{n=0}^{N-1}\sum_{n'=0}^{N-1}\cov\left\lbrace \bfW_n(\nu),\bfW_{n'}(\nu)\right\rbrace \nonumber\\
&=\lim_{N\rightarrow\infty}\frac{1}{N}\sum_{h=-N+1}^{N-1}\left( N-|h|\right)\cov\left\lbrace \bfW_n(\nu),\bfW_{n-h}(\nu)\right\rbrace=\sum_{h=-\infty}^{\infty}\Lambda_h(\nu).
\end{align}
Under the assumptions of Theorem \ref{thm}, the moments $\Lambda_h(\nu)$ are consistently estimated by $\hat{\Lambda}_h(\nu)$, for $0\leq h< N$,
$$\hat{\Lambda}_h(\nu)=\frac{1}{N}\sum_{n=0}^{N-h-1}\hat{\bfW}_{n}(\nu)\hat{\bfW}_{n-h}^\top(\nu)\quad\text{and}\quad\hat{\Lambda}_{-h}(\nu)=\hat{\Lambda}^{\top}_h(\nu).$$
This raises the question of whether matrix $$\check{\bfPsi}(\nu)=\sum_{h=-N+1}^{N-1}\hat{\Lambda}_h(\nu)$$ would be a consistent estimator of $\bfPsi(\nu)$. The answer is clearly negative since, for all $N$,
$$\check{\bfPsi}(\nu)=\frac{1}{N}\left(\sum_{n=0}^{N-1}\hat{\bfW}_{n}(\nu)
\right)^2=\frac{1}{N}\left(\sum_{n=0}^{N-1}\bfX_n(\nu)\otimes\hat{\bfepsilon}_{ns+\nu}\right)^2=\frac{1}{N}\left(\frac{\partial S(\hat{\bfbeta})}{\partial\bfbeta}\right)^2=0.$$
Note that when the index $|h|$ in \eqref{I-noyau} is large, the moments $\Lambda_h(\nu)$ are likely to be poorly estimated since their estimators are based on only few observations. The classical solution to get around this problem is to weight the empirical moments $\hat{\Lambda}_{h}(\nu)$.
To estimate $\bfPsi(\nu)$, we consider a sequence of real numbers $(b_N)_{N\in\mathbb{N}}$ such that
\begin{equation}\label{conditions_bn}b_N\to 0\quad\text{and}\quad
Nb_N^{\frac{10+4\kappa}{\kappa}}\to \infty\quad\text{as}\quad
N\to\infty,
\end{equation} and a weight function
$f:\mathbb{R}\to\mathbb{R}$ which is bounded, with compact support
$[-a,a]$ and continuous at the origin with $f(0)=1$. Note that under
the above assumptions, we have
\begin{equation}\label{fenetreO(1)}
b_N\sum_{|h|<N}|f(hb_N)|=\mathrm{O}(1).
\end{equation}
% by weights of the form $f(hb_N)$, where $(b_N)_{N\geq 0}$ is a sequence of real numbers and $f: \mathbb{R}\to \mathbb{R}$ is assumed to be a bounded function, with compact support $[-a,a]$ and continuous at the origin with $f(0)=1$.
Examples of such weight functions can be found in \cite{BM14} and are:  the truncated uniform or rectangular (REC) window $f(x)=\mathds{1}_{[-1,1]}(x)$, the Bartlett (BAR) window $f(x)=(1-|x|)\mathds{1}_{[-1,1]}(x)$, the quadratic-spectral window $$f(x) =\dfrac{25}{12\pi^2x^2}\left(\dfrac{\sin{(6\pi x/5)}}{6\pi x/5}-\cos{(6\pi x/5)}\right) $$ or the Parzen (PAR) window
\begin{equation}\label{Parzen}
f(x)={\begin{cases}
  1-6x^2+6|x|^3 & \quad \text{ if } |x|\leq 1/2\\
  2(1-|x|)^3 & \quad \text{ if } 1/2\leq |x|\leq 1\\
 	0 & \quad \text{ otherwise }\\ \end{cases}} \ .
\end{equation}
Consider the matrix
$$\hat{\bfPsi}^\mathrm{HAC}(\nu):=\sum_{h=-T_N}^{T_N}f(hb_N)\hat{\Lambda}_{h}(\nu) \quad\text{and}\quad T_N=\left\lfloor\frac{a}{b_N}\right\rfloor,$$
where $\lfloor x \rfloor$  denotes the integer part of the real $x$.

\noindent We are now able to state the following theorem, which shows
the weak consistency of an empirical estimator of  $\hat{\bfPsi}^\mathrm{HAC}(\nu)$.
\begin{theorem}\label{estHAC}
Under the assumptions of Theorem \ref{thm} and  if the sequence $(b_N)_{N\geq 0}$ is chosen such that \eqref{conditions_bn} is satisfied, we have
\begin{eqnarray*}
\hat{\bfPsi}^\mathrm{HAC}(\nu)\cp \bfPsi(\nu)\text{  as } N\to\infty.
\end{eqnarray*}
\end{theorem}
\noindent The proof of Theorem \ref{estHAC} is postponed to Section \ref{app}.

\noindent Theorems~\ref{convergence_Isp}  and \ref{estHAC}, and \eqref{estim_J} show that
\begin{align}
\label{estThetaSP}
  \hat{\bfTheta}^{\mathrm{SP}}(\nu)&:=\left(\hat{\bfOmega}^{-1}(\nu) \otimes \bfI_d\right)\hat{\bfPsi}^{\mathrm{SP}}(\nu)\left(\hat{\bfOmega}^{-1}(\nu) \otimes \bfI_d\right)
  \\ \text{and}\quad
    \hat{\bfTheta}^{\mathrm{HAC}}(\nu)&:=\left(\hat{\bfOmega}^{-1}(\nu) \otimes \bfI_d\right)\hat{\bfPsi}^{\mathrm{HAC}}(\nu)\left(\hat{\bfOmega}^{-1}(\nu) \otimes \bfI_d\right)\label{estThetaHAC}
\end{align}
\noindent are weakly consistent estimators of $\bfTheta(\nu)$.

\section{Testing linear restrictions about the parameter}
\label{test}
\noindent In addition to the $K(\nu)$ linear constraints imposed in Section \ref{LScont},
the parameter may satisfy other linear constraints  which can be interesting to test (in particular $\bfPhi_{p(\nu)}=0$, $\nu = 1,\ldots,s$). Theorems \ref{convergence_Isp}-\ref{estHAC} and \eqref{NAcont} can be exploited to test $s_0(\nu)$ linear constraints on the elements of the free parameter $\bfxi(\nu)$. The null hypothesis takes the form $$H_0:\bfR_0(\nu)\bfxi(\nu)=\bfr_0(\nu), \qquad\nu = 1,\ldots,s$$
where $\bfR(\nu)$ is a known $\{ s_0(\nu) \} \times K(\nu)$
matrix of rank $s_0(\nu)$ and $\bfr_0(\nu)$ is a known $\{ s_0(\nu) \} $-dimensional vector. The Wald  principle is employed frequently for testing $H_0$. We now examine if this principle
remains valid in the non standard framework of weak PVAR models.

\noindent From \eqref{NAcont}, we deduce that
\begin{align}
\label{NAcont2}
N^{1/2} \{  \bfR_0(\nu)\hat{\bfxi}(\nu) - \bfr_0(\nu) \} \cd
  N_{s_0(\nu)}\left(\bfzero,
              \bfR_0(\nu)\bfTheta^{\bfxi}(\nu)\bfR_0^\top(\nu)
            \right).
\end{align}
Let
\begin{align}\nonumber
  \hat{\bfTheta}^{\bfxi}(\nu)&=\left[ \bfR^{\top}(\nu) \{ \hat{\bfOmega}(\nu) \otimes \hat{\bfSigma}^{-1}_{\bfepsilon}(\nu) \} \bfR(\nu) \right]^{-1}
\bfR^{\top}(\nu) \{ \bfI_{dp(\nu)} \otimes \hat{\bfSigma}^{-1}_{\bfepsilon}(\nu) \}  \hat{\bfPsi}(\nu) \\&\qquad\times\{ \bfI_{dp(\nu)} \otimes \hat{\bfSigma}^{-1}_{\bfepsilon}(\nu) \}   \bfR(\nu)
 \left( \left[ \bfR^{\top}(\nu) \{ \hat{\bfOmega}(\nu) \otimes \hat{\bfSigma}^{-1}_{\bfepsilon}(\nu) \} \bfR(\nu) \right]^{-1}\right)^{\top}\label{estThetacont}
\end{align}
be a consistent estimator of $\bfTheta^{\bfxi}(\nu)$, where $\hat{\bfOmega}(\nu)$, $\hat{\bfSigma}_{\bfepsilon}(\nu)$ and $\hat{\bfPsi}(\nu)$
are consistent estimators of $\bfOmega(\nu)$, $\bfSigma_{\bfepsilon}(\nu)$ and $\bfPsi(\nu)$, as defined in Section \ref{estimOmega}.
In view of \eqref{NAcont2}, under the assumptions of Theorems~\ref{thm}, \ref{convergence_Isp} and \ref{estHAC}, and the assumption that $\bfPsi(\nu)$ is invertible, the modified Wald statistic
\begin{equation}
\label{modwald}
\mathrm{{W}}_N(\nu)=N\left(  \bfR_0(\nu)\hat{\hat{\bfxi}}(\nu) - \bfr_0(\nu) \right)^\top\left(\bfR_0(\nu)\hat{\bfTheta}^{\bfxi}(\nu)\bfR_0^\top(\nu)\right)^{-1}\left(  \bfR_0(\nu)\hat{\hat{\bfxi}}(\nu) - \bfr_0(\nu) \right)
\end{equation}
asymptotically follows a $\chi^2_{s_0(\nu)}$ distribution under $H_0$. Therefore, the standard formulation of the Wald test remains valid.
More precisely, at the asymptotic level $\alpha$,
the modified Wald test consists in rejecting $H_0$  when
$\mathrm{{W}}_N(\nu)>\chi^2_{s_0(\nu)}(1-\alpha)$. It is however important to note that a consistent estimator of the form \eqref{estThetacont} is required.
Note that in the strong PVAR case and in view of Remark \ref{rem2}, $\hat{\bfTheta}^{\bfxi}(\nu)=\hat{\bfTheta}_S^{\bfxi}(\nu)=:\left[ \bfR^{\top}(\nu) \{ \hat{\bfOmega}(\nu) \otimes \hat{\bfSigma}^{-1}_{\bfepsilon}(\nu) \} \bfR(\nu) \right]^{-1}$ and the Wald statistic takes the more conventional form:
\begin{equation}
\label{standwald}
\mathrm{{W}}^\ast_N(\nu)=:N\left(  \bfR_0(\nu)\hat{\hat{\bfxi}}(\nu) - \bfr_0(\nu) \right)^\top\left(\bfR_0(\nu)\hat{\bfTheta}_S^{\bfxi}(\nu)\bfR_0^\top(\nu)\right)^{-1}\left(  \bfR_0(\nu)\hat{\hat{\bfxi}}(\nu) - \bfr_0(\nu) \right)
\end{equation}
The estimator $\hat{\bfTheta}_S^{\bfxi}(\nu)$ of $\bfTheta_S^{\bfxi}(\nu)$, which is routinely used in the time series software, is only valid in the strong PVAR case.

%%% ----------------------------------------------------------------------
%%% ----------------------------------------------------------------------
\section{Simulations}\label{simul}
\noindent By means of a small Monte Carlo experiment, we investigate the behaviour of the least squares estimators for strong and weak bivariate PVAR model. The following data generating process (DGP) is used:
\begin{eqnarray}\label{dgp}
\mbox{DGP} &:& \bfY_{ns+\nu} = \bfPhi(\nu)\bfY_{ns+\nu-1} + \bfepsilon_{ns+\nu},\quad \nu = 1,\ldots,5.
\end{eqnarray}
We considered the case of five seasons, that is $s=5$. The model DGP corresponds to a PVAR model of order $1$. The coefficients of the DGP  in \eqref{dgp} are chosen such that Assumption \textbf{(A1)}  holds
and are given in  Table~\ref{DGP1}.

\begin{table}[H]
\caption{Parameters of DGP models used in the simulation}\label{DGP1}
\vspace{0.1cm}
\centering
\begin{tabular}{l cc cccc cc cccc cc cccc cc cccc cc cccc}
  \toprule
  MODEL & & \multicolumn{2}{c}{$\bfPhi(1)$} & & \multicolumn{2}{c}{$\bfPhi(2)$}& & \multicolumn{2}{c}{$\bfPhi(3)$}& & \multicolumn{2}{c}{$\bfPhi(4)$}& & \multicolumn{2}{c}{$\bfPhi(5)$}\\
  \midrule
  DGP & &-1.43& 0.00&& 0.46& 0.00&&1.23& 0.00&&0.30& 0.00&&0.90& 0.00\\
          & &0.00& 0.62&& 0.00& 0.70&&0.00& -0.30&&0.00& 0.45&&0.00& 0.20\\
  \bottomrule
\end{tabular}
\end{table}

\noindent We consider that  the stochastic process $\bfepsilon = \{ \bfepsilon_t, t \in \mathbb{Z} \}$ in \eqref{dgp} corresponds to a zero mean periodic white noise with the error covariance matrix $\bfSigma_{\bfepsilon}(\nu)$ given in Table~\ref{error}.

\begin{table}[H]
\caption{Error covariance matrices used in the simulation}\label{error}
\vspace{0.1cm}
\centering
\begin{tabular}{c cccc cc cccc cc cccc cc cccc cc cccc}
  \toprule
    & \multicolumn{2}{c}{$\bfSigma_{\bfepsilon}(1)$} & & \multicolumn{2}{c}{$\bfSigma_{\bfepsilon}(2)$}& & \multicolumn{2}{c}{$\bfSigma_{\bfepsilon}(3)$}& & \multicolumn{2}{c}{$\bfSigma_{\bfepsilon}(4)$}& & \multicolumn{2}{c}{$\bfSigma_{\bfepsilon}(5)$} \\
  \midrule
           & 1.00 & 0.05 & & 1.60 & 0.30&& 2.20& -0.20&&2.50& -0.10&&0.90& 0.00 \\
           & 0.05 & 1.50 & & 0.30 & 0.50&&-0.20& 0.80&&-0.10& 1.20&& 0.00& 1.70 \\
  \bottomrule
\end{tabular}
\end{table}
\noindent First we study numerically the behavior of the least squares estimators for strong and weak PVAR models of the form \eqref{dgp}.
We consider the strong PVAR case by assuming that the innovation process $\bfepsilon$ in \eqref{dgp} is defined by an iid sequence such that
\begin{equation} \label{bruitfort}
\bfepsilon_{ns+\nu} =\left(\begin{array}{c}\epsilon_{1,ns+\nu}\\\epsilon_{2,ns+\nu}\end{array}\right) \  \overset{\text{law}}{=} \ {\cal
N}(0,I_2).
\end{equation}

\noindent We repeat the same experiment on a weak PVAR model, meaning that the stochastic process $\bfepsilon$  defined by \eqref{weaknoise} and given for $m=2$ by
\begin{equation}
\label{bruitfaible}
\bfepsilon_{ns+\nu} =\bfM_{\nu}^\top
\begin{pmatrix}
\eta_{1,ns+\nu}\eta_{1,ns+\nu-1}\eta_{1,ns+\nu-2} \\
\eta_{2,ns+\nu}\eta_{2,ns+\nu-1}\eta_{2,ns+\nu-2}
\end{pmatrix}
\quad\text{where}\quad\bfM^{T}_{\nu} \bfM_{\nu} = \bfSigma_{\bfepsilon}(\nu).
\end{equation}

\begin{figure}[H]
\begin{subfigure}{.5\textwidth}
  \centering
  % include first image
  \includegraphics[width=1\linewidth]{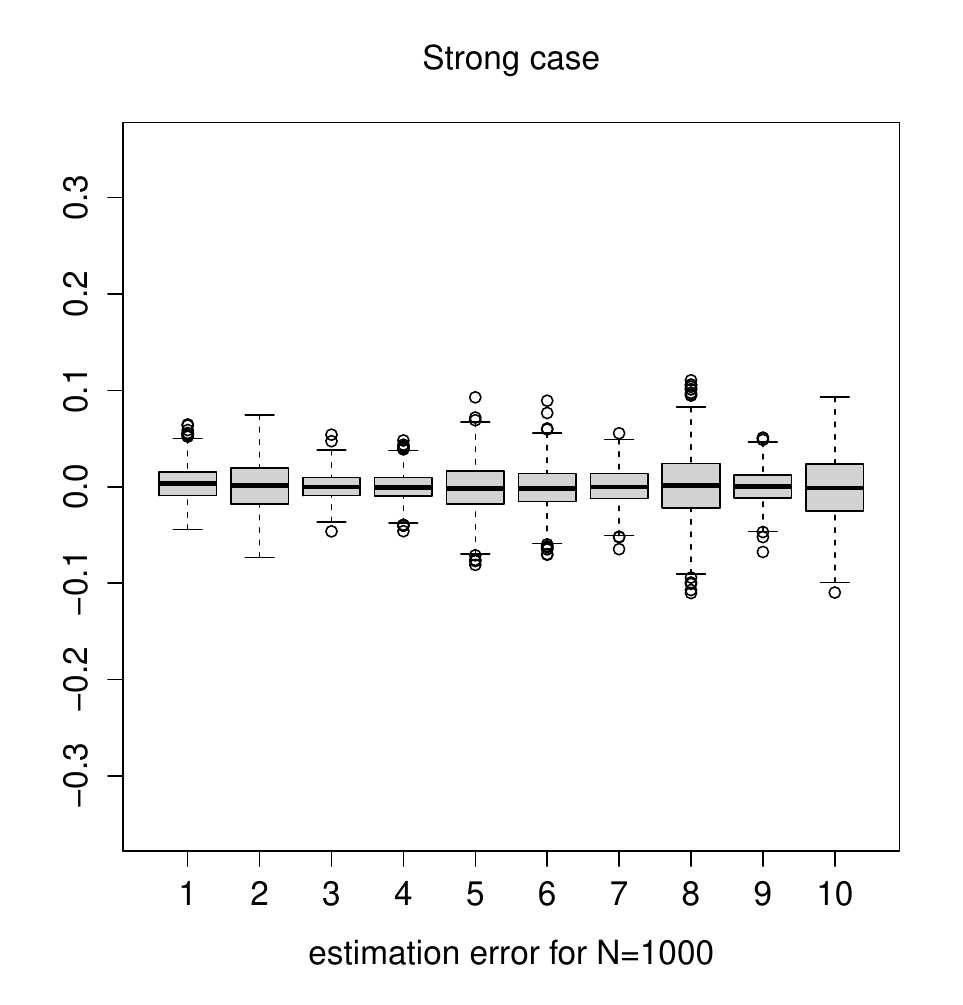}
  %\caption{Estimates of diag($\Theta_S(\nu)$) (Remark~3.2)}
  %\label{fig:sub-first}
\end{subfigure}
\begin{subfigure}{.5\textwidth}
  \centering
  % include second image
  \includegraphics[width=1\linewidth]{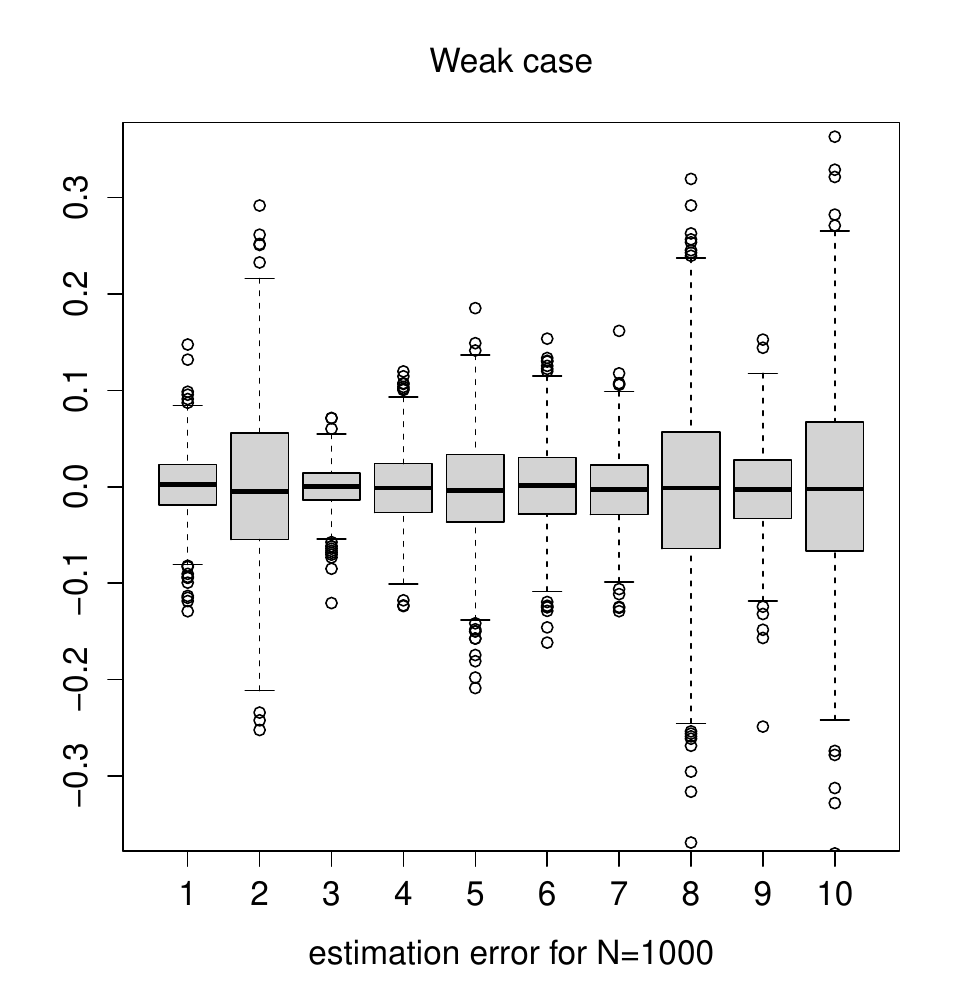}
  %\caption{Estimates of diag($\Theta^{SP}(\nu)$) (Equation~42)}
  %\label{fig:sub-second}
\end{subfigure}
\caption{The least squares estimator of $1000$ independent simulations of model \eqref{dgp} with unknown parameter given in Table~\ref{DGP1}, when the noise is strong (left panels) and when the noise is weak (right panels).
The panels display the distribution of the estimation errors $\hat{\bfPhi}_{ii}(\nu)-\bfPhi_{ii}(\nu)$, for $i=1,2$ and $\nu=1,\dots,5$.}
\label{fig:errors}
\end{figure}

\begin{figure}[H]
\begin{subfigure}{.5\textwidth}
  \centering
  % include first image
  \includegraphics[width=1\linewidth]{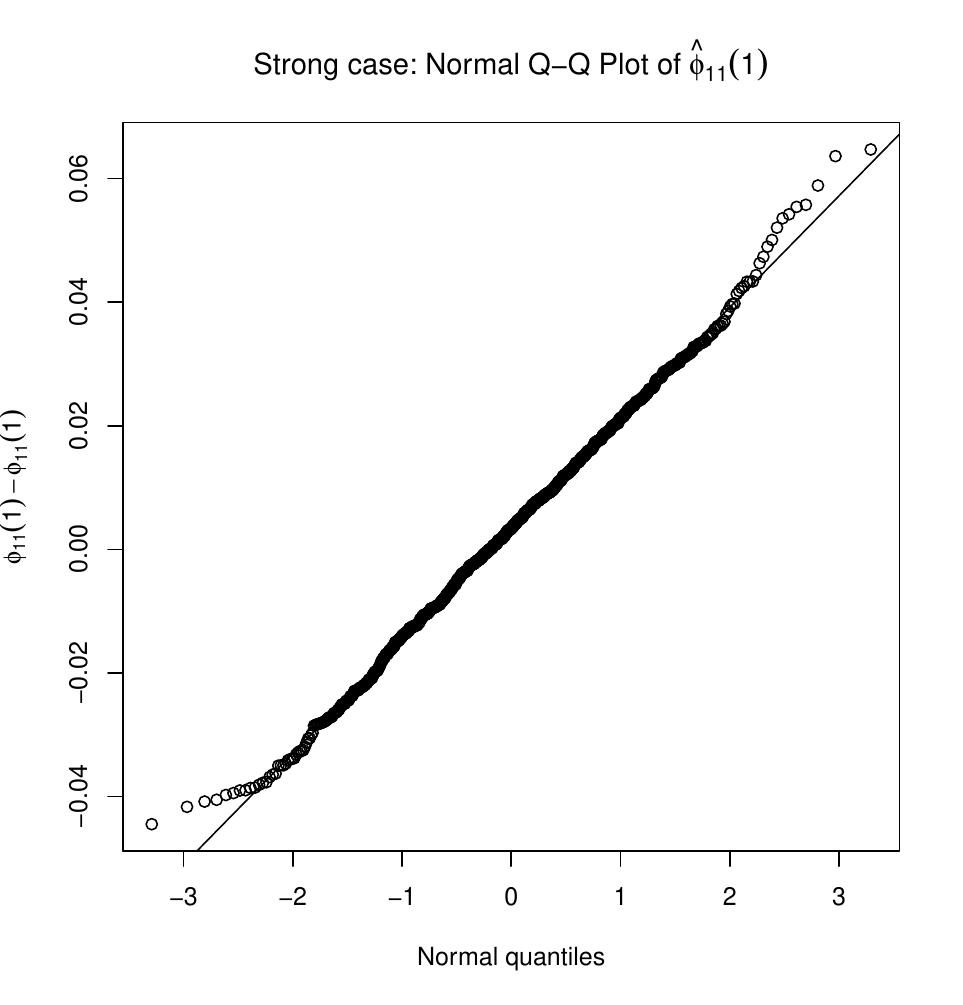}
  %\caption{Estimates of diag($\Theta_S(\nu)$) (Remark~3.2)}
  %\label{fig:sub-first}
\end{subfigure}
\begin{subfigure}{.5\textwidth}
  \centering
  % include second image
  \includegraphics[width=1\linewidth]{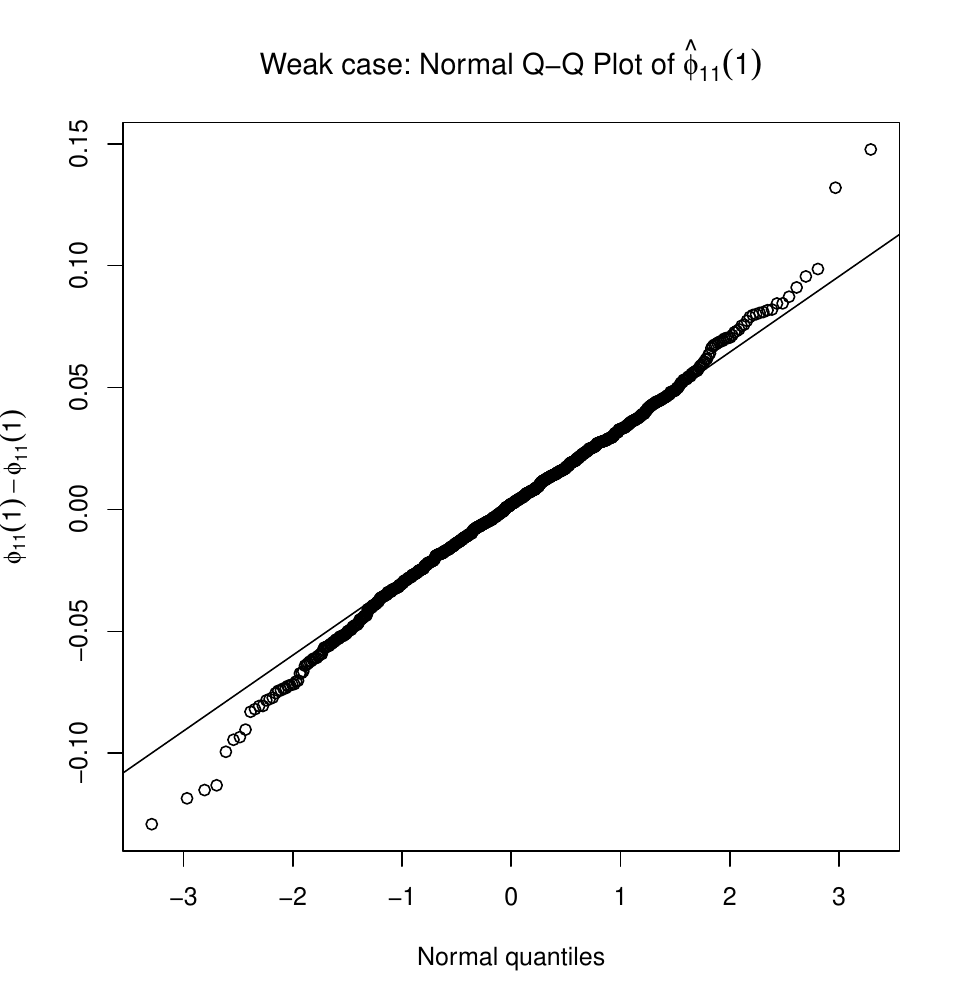}
  %\caption{Estimates of diag($\Theta^{SP}(\nu)$) (Equation~42)}
  %\label{fig:sub-second}
\end{subfigure}

\begin{subfigure}[H]{.5\textwidth}
  \centering
  % include third image
  \includegraphics[width=1\linewidth]{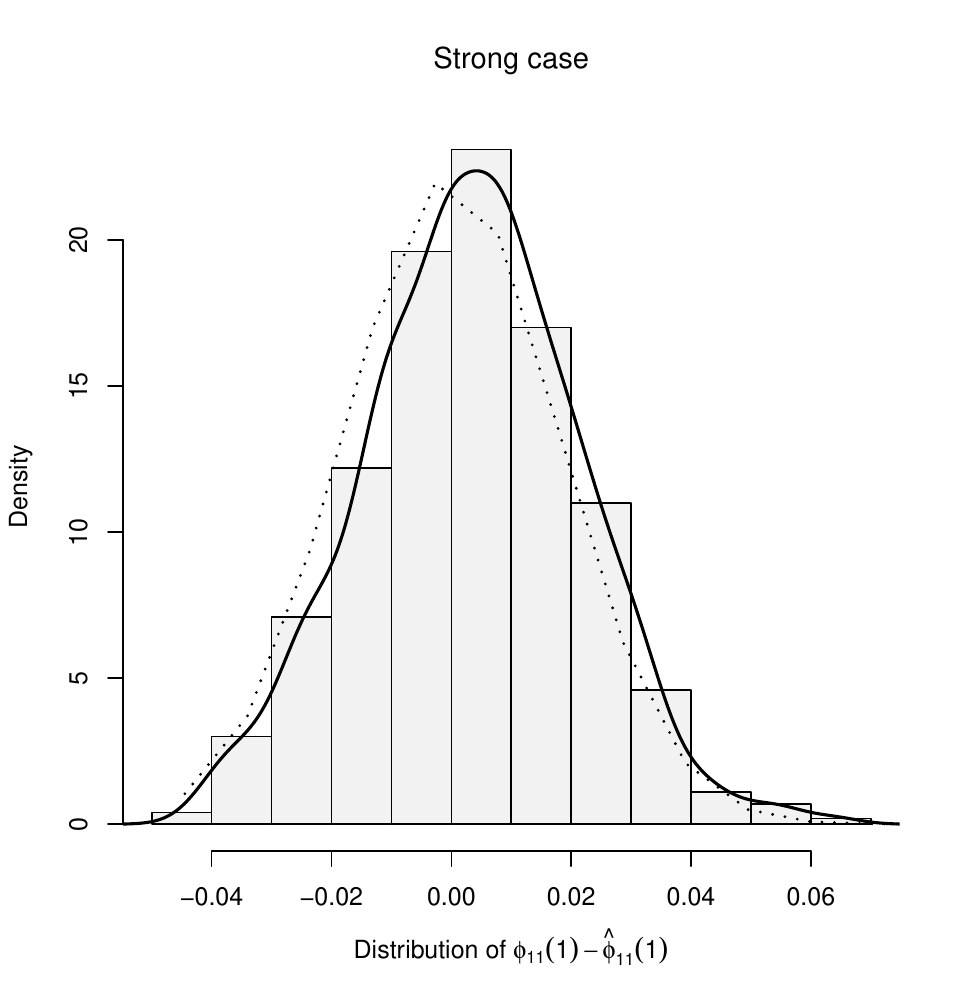}
  %\caption{Estimates of diag($\Theta_S(\nu)$) (Remark~3.2)}
  %\label{fig:sub-third}
\end{subfigure}
\begin{subfigure}{.5\textwidth}
  \centering
  % include fourth image
  \includegraphics[width=1\linewidth]{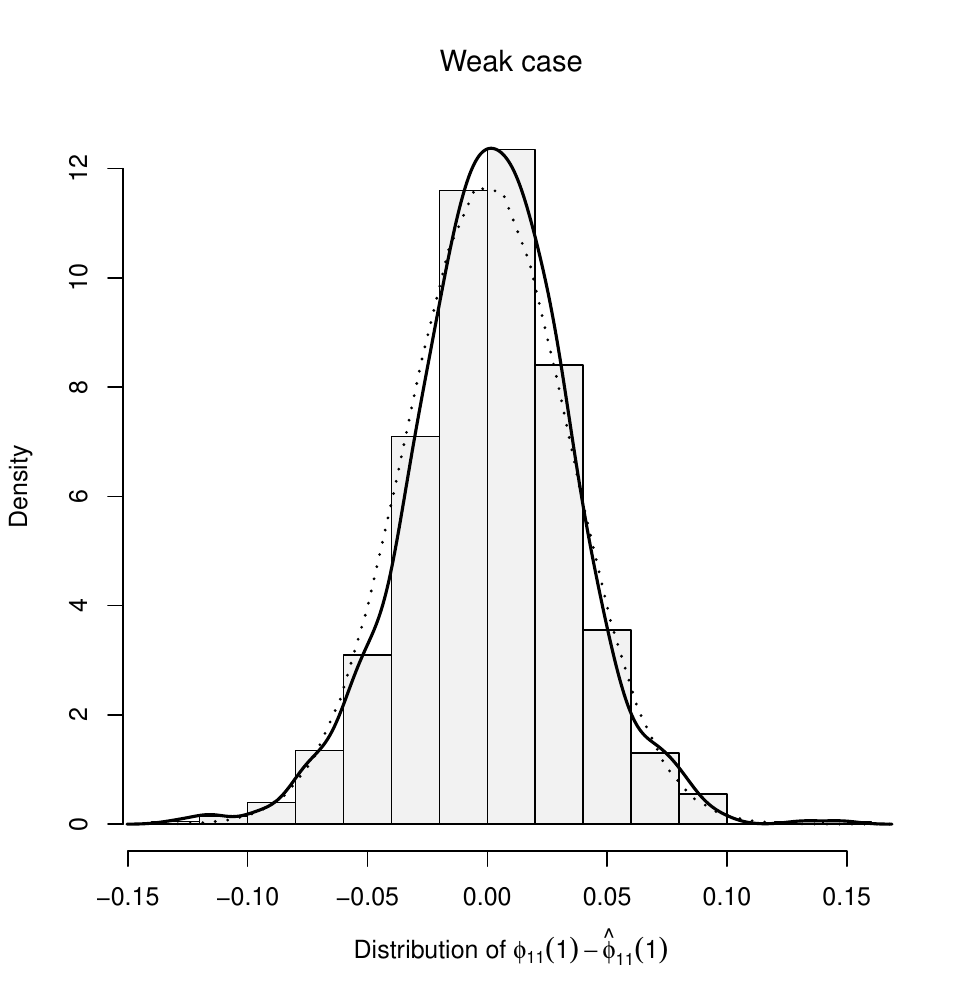}
  %\caption{Estimates of diag($\Theta^{SP}(\nu)$) (Equation~42)}
  %\label{fig:sub-fourth}
\end{subfigure}
\caption{The least squares estimator of $1000$ independent simulations of the model \eqref{dgp} with size $N=1000$ and unknown parameter given in Table~\ref{DGP1}, when the noise is strong (left panels) and when the noise is weak (right panels).
The panels of the top present the Q–Q plot of the estimates $\bfPhi_{11}(1)$. The bottom panels display the distribution of the same estimates. The kernel density estimate is displayed in full line, and the centered Gaussian density with the same variance is plotted in dotted line.
}
\label{fig:estimates}
\end{figure}

\noindent Figures~\ref{fig:errors} and~\ref{fig:estimates} compare the distribution of the least squares estimator in the strong and weak noise cases. The distributions of $\hat{\bfPhi}_{ii}(\nu)$, for $i=1,2$ and $\nu=1,\dots,5$ are more accurate in the strong case than in the weak one. Similar simulation experiments, not reported here, reveal that the situation is opposite, that is the least squares estimators of $\hat{\Phi}_{ii}(\nu)$ are more accurate in the weak case than in the strong case, when the weak noise is defined by
$\epsilon_{i,ns+\nu}=\eta_{i,ns+\nu}(|\eta_{i,ns+\nu-1}|+1)^{-1}$ for $i=1,2$. This is in accordance with the univariate results of \cite{RT96} who showed that, with similar noises, the asymptotic variance of the sample autocorrelations can be greater or less than 1 as well (1 is the asymptotic variance for strong white noises).
Figure~\ref{fig:estimates} compares the distribution of  $\bfPhi_{11}(1)$ in the strong and weak noise cases. We consider here (one of) the parameter which variance's seems to have problems in the weak case, when we use the standard estimator $\hat{\bfTheta}_{S}(\nu)$.

\begin{sidewaysfigure}
\begin{subfigure}{.33\textwidth}
  \centering
  % include first image
  \includegraphics[width=1.05\linewidth]{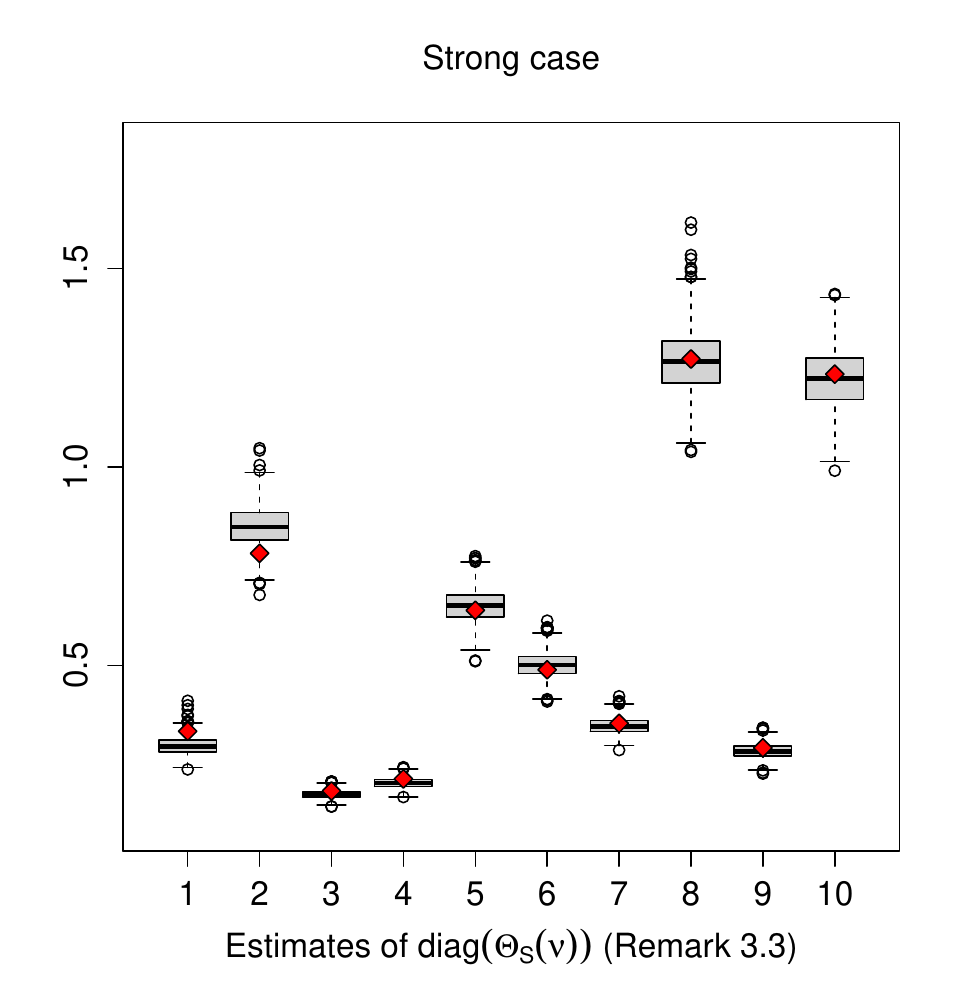}
  %\caption{Estimates of diag($\Theta_S(\nu)$) (Remark~3.2)}
  %\label{fig:sub-first}
\end{subfigure}
\begin{subfigure}{.33\textwidth}
  \centering
  % include fourth image
  \includegraphics[width=1.05\linewidth]{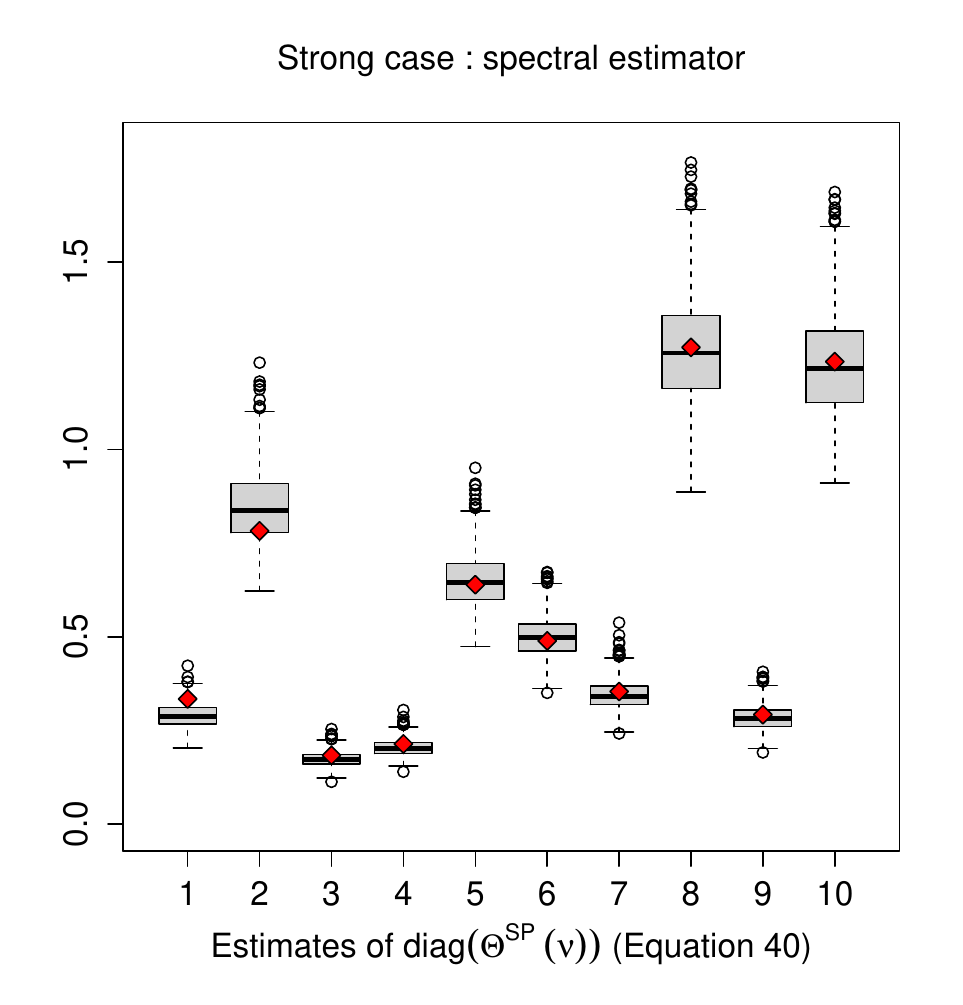}
  %\caption{Estimates of diag($\Theta^{SP}(\nu)$) (Equation~42)}
  %\label{fig:sub-fourth}
\end{subfigure}
\begin{subfigure}{.33\textwidth}
  \centering
  % include second image
  \includegraphics[width=1.05\linewidth]{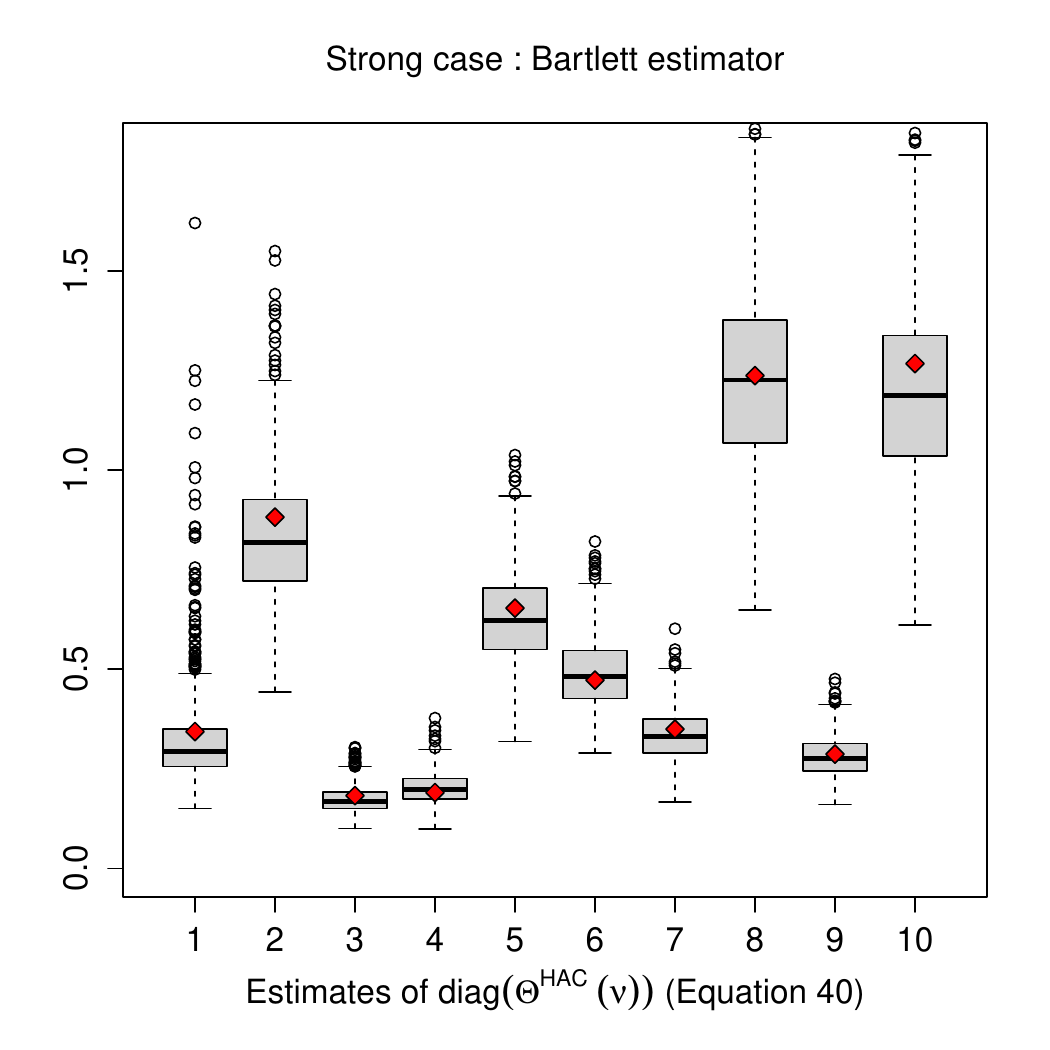}
  %\caption{Estimates of diag($\Theta^{SP}(\nu)$) (Equation~42)}
  %\label{fig:sub-second}
\end{subfigure}\hspace{1em}

\begin{subfigure}{.33\textwidth}
  \centering
  % include third image
  \includegraphics[width=1.05\linewidth]{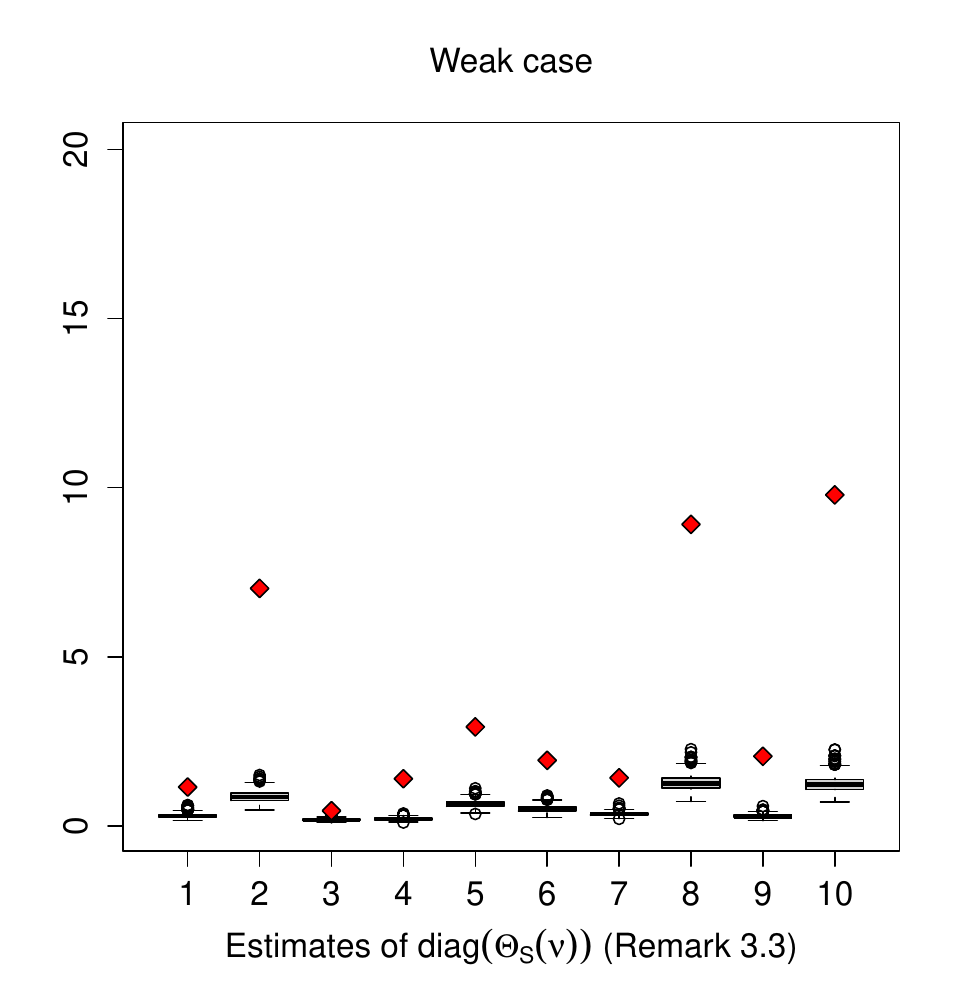}
  %\caption{Estimates of diag($\Theta_S(\nu)$) (Remark~3.2)}
  %\label{fig:sub-third}
\end{subfigure}
\begin{subfigure}{.33\textwidth}
  \centering
  % include second image
  \includegraphics[width=1.05\linewidth]{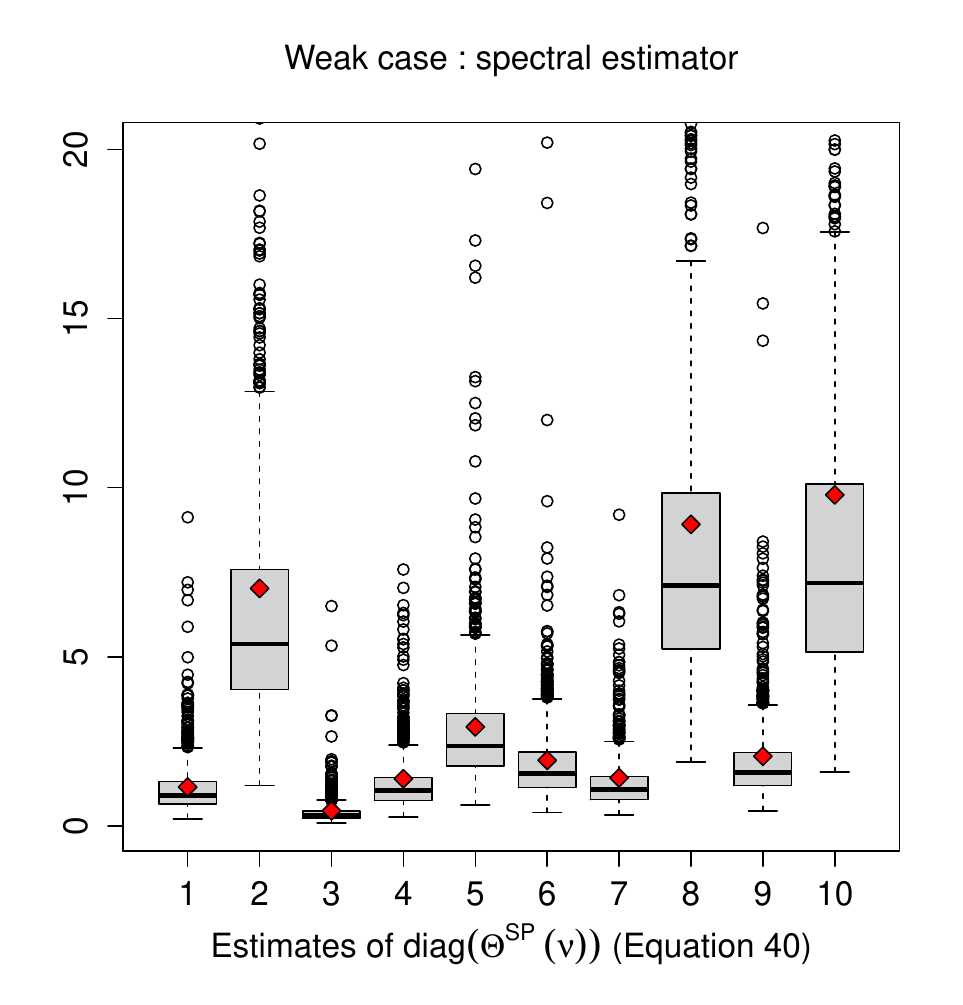}
  %\caption{Estimates of diag($\Theta^{SP}(\nu)$) (Equation~42)}
  %\label{fig:sub-second}
\end{subfigure}
\begin{subfigure}{.33\textwidth}
  \centering
  % include second image
  \includegraphics[width=1.05\linewidth]{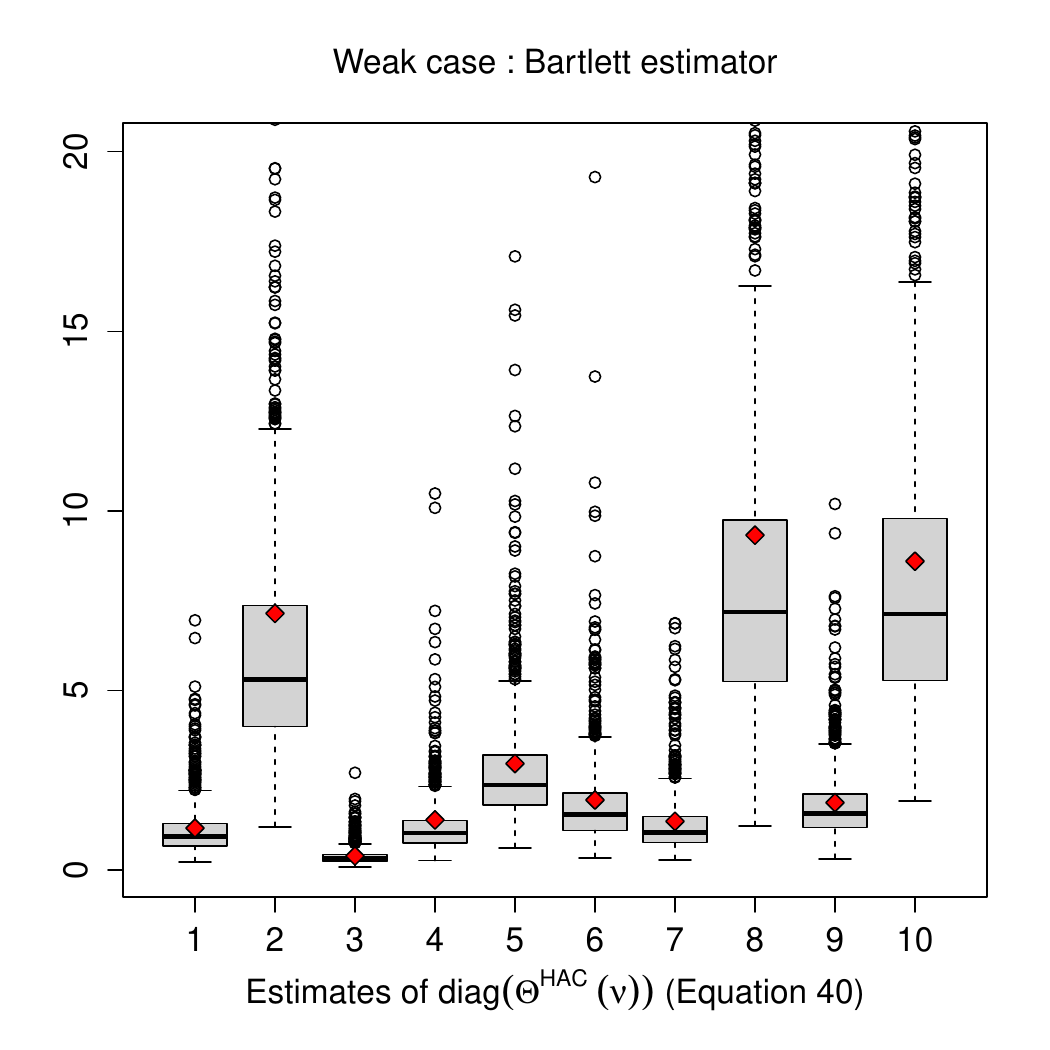}
  %\caption{Estimates of diag($\Theta^{SP}(\nu)$) (Equation~42)}
  %\label{fig:sub-second}
\end{subfigure}
\caption{Comparison of standard and modified estimates of the asymptotic variance of the least squares estimator of the model parameters, on the simulated models presented in~Figure~\ref{fig:errors}. The diamond symbols represent the mean, over $1000$ replications, of the standardized squared errors $N\left(\hat{\bfPhi}_{11}(1)-\bfPhi_{11}(1)\right)^2$ for 1 (0.33 in the strong
case and 1.15 in the weak case), $N\left(\hat{\bfPhi}_{22}(1)-\bfPhi_{22}(1)\right)^2$ for 2 (0.78 in the strong
case and 7.02 in the weak case),
$N\left(\hat{\bfPhi}_{11}(2)-\bfPhi_{11}(2)\right)^2$ for 3 (0.18 in the strong
case and 0.45 in the weak case), $N\left(\hat{\bfPhi}_{22}(2)-\bfPhi_{22}(2)\right)^2$ for 4 (0.21 in the strong
case and 1.39 in the weak case),
$N\left(\hat{\bfPhi}_{11}(3)-\bfPhi_{11}(3)\right)^2$ for 5 (0.64 in the strong
case and 2.93 in the weak case), $N\left(\hat{\bfPhi}_{22}(3)-\bfPhi_{22}(3)\right)^2$ for 6 (0.49 in the strong
case and 1.94 in the weak case),
$N\left(\hat{\bfPhi}_{11}(4)-\bfPhi_{11}(4)\right)^2$ for 7 (0.35 in the strong
case and 1.42 in the weak case), $N\left(\hat{\bfPhi}_{22}(4)-\bfPhi_{22}(4)\right)^2$ for 8 (1.27 in the strong
case and 8.91 in the weak case),
$N\left(\hat{\bfPhi}_{11}(5)-\bfPhi_{11}(5)\right)^2$ for 9 (0.29 in the strong
case and 2.06 in the weak case) and $N\left(\hat{\bfPhi}_{22}(5)-\bfPhi_{22}(5)\right)^2$ for 10 (1.23 in the strong
case and 9.79 in the weak case) .}
\label{fig:var}
\end{sidewaysfigure}

\noindent Figure~\ref{fig:var} compares the standard estimator $\hat{\bf\Theta}_S(\nu)$ with the proposed sandwich estimators based on spectral density estimation $\hat{\bf\Theta}^\mathrm{{SP}}(\nu)$ or on kernel methods $\hat{\bf\Theta}^\mathrm{{HAC}}(\nu)$ of the asymptotic variance $\bf\Theta(\nu)$. We used the spectral estimator $\hat{\bfPsi} = \hat{\bfPsi}^\mathrm{{SP}}$ defined in Theorem~\ref{convergence_Isp} where the AR order $r$ is automatically selected by AIC, using the function \verb"VARselect()" of the \textit{vars} R package. Note that similar simulation experiments, not reported here, reveal that the performance of the proposed estimator is least sensitive to the choice of others criteria  such that: BIC, HQ and FPE. The HAC estimator based on kernel methods $\hat{\bfPsi} = \hat{\bfPsi}^\mathrm{{HAC}}$ defined in Theorem~\ref{estHAC} is also used. HAC estimators have been the focus of extensive research in the time series literature. Contributions to this research in the econometrics literature include, among others,~\cite{newey},~\cite{A91},~\cite{Mu14}  and~\cite{LLSW18}.The bandwidth selection for the HAC estimation is an important practical issue.
For kernel densities with unit-interval support, the bandwidth parameter, is often called the lag-truncation parameter. Based on theoretical results in~\cite{A91}, the practice is to choose a small value for the lag-truncation parameter. More recently, it has been shown that this standard approach can often lead to tests which incorrectly reject the null hypothesis~\citep{Mu14}. Much of the literature remains in the Newey-West framework but uses very long lag-truncation parameter~\citep{KV02}. As indicated by~\cite{FZ07}, "it is well known that choice of bandwidth equal to the sample size (i.e. $b_N=1/N$ in our case) results in inconsistent LRV estimators". In our case it is crucial to have a consistent estimator of the matrix $\bfPsi (\nu)$.

Several leading lag-truncation choices based on traditional Newey-West HAC estimators are:
\begin{itemize}
\item $b_N = 1/\ln{(N)}$, as proposed by~\cite{FZ00}. 
%For example, if $N = 1000$, then $b_N = 1/6$.
\item $b_N = 1/\left(\lfloor 4(N/100)^{2/9}\rfloor+1\right)$ or $b_N = 1/\left(\lfloor N^{1/4}\rfloor+1\right)$. 

This choice is a standard textbook recommendation~\citep{Wo15}.
\item $b_N = 1/\left(\lfloor 0.75N^{1/3}\rfloor+1\right)$. This rule derives from a formula of~\cite{A91}, in the case of a first order autoregressive model.
\item $ b_N=1/\left(\lfloor 1.3N^{1/2}\rfloor + 1\right)$, as proposed by~\cite{LLSW18}. Its use of $N^{1/2}$  produces higher truncation lags. For example, if $N = 1000$, then $b_N = 1/43$.
\item $b_N = 1/N$, as proposed by~\cite{KV02}.
\end{itemize}
The performance of Newey-West estimators depends on the choice of the kernel function and lag-truncation. We focused our investigation to weight function
$f:\mathbb{R}\to\mathbb{R}$ which is bounded, with compact support
$[-a,a]$ and continuous at the origin with $f(0)=1$. Such weight functions are for instance: Bartlett, Truncated, Parzen and Quadratic Spectral kernels. To determine the optimal lag-truncation parameter and the kernel, a $10$-fold cross-validation is used. In the density estimation literature, the cross-validation method has been suggested by~\cite{BB87} or by~\cite{WW75}. For the bandwidth, we chose $30$ values between $1/50$ and $1/6$ and four kernels were investigated: Bartlett, Truncated, Parzen and Quadratic Spectral. The best results, for our simulated data, was obtained using the Bartlett kernel with a bandwidth equal to $1/21$.

In the strong PVAR case we know that the three estimators are consistent. In view of the three top panels of Figure~\ref{fig:var}, it seems that the standard estimator is most accurate than the proposed sandwich  estimators in the strong case. This is not surprising because the spectral estimator or the HAC estimator are more robust, in the sense that these estimators continue to be consistent in the weak PVAR case, contrary to the standard estimator. It is clear that in the weak case $N Var(\hat{\bfPhi}_{ii}(\nu)-\bfPhi_{ii}(\nu))^2$ is better estimated by $\hat{\bfTheta}_{ii}^\mathrm{{SP}}(\nu)$ or by $\hat{\bf\Theta}_{ii}^\mathrm{{HAC}}(\nu)$ (see the box-plots $1,2,\dots,10$ of the center-bottom and the right-bottom panel of Figure~\ref{fig:var}) than by $\hat{\bfTheta}_{S}(\nu)$ (see the box-plots $1,2,\dots,10$ of the left-bottom panel), for $i=1,2$ and $\nu=1,\dots,5$. The failure of the standard estimator of $\bfTheta$ in the weak PVAR setting may have important consequences in terms of hypothesis testing for instance.

\noindent Table~\ref{tabwald} displays the empirical sizes  of the standard Wald test and that of the modified versions proposed
in Section \ref{test}. We use  $3$ nominal levels $\alpha=1\%$, 5\% and 10\%. For these nominal levels, the empirical  relative frequency of rejection size over the $1000$ independent replications should vary  respectively within the confidence intervals $[0.3\%,1.7\%]$, $[3.6\%, 6.4\%]$  and $[8.1\%,11.9\%]$ with probability 95\% and  $[0.3\%, 1.9\%]$, $[3.3\%, 6.9\%]$ and $[7.6\%, 12.5\%]$ with probability 99\% under the assumption that the true probabilities of rejection are respectively $\alpha=1\%$, $\alpha=5\%$ and $\alpha=10\%$.
When the relative rejection frequencies  are outside the
significant limits with probability 95\%, they are displayed in bold type in Table \ref{tabwald}. For the strong PVAR model I, the relative rejection frequencies are inside the significant limits. For the weak PVAR model II, the relative rejection frequencies of the standard Wald  test are definitely outside the significant limits. It may lead the statistician to wrongly reject the hypothesis that $H_0:\bfPhi_{22}(\nu)=0$ for $\nu=1,\dots,5$ if he does not take into account the dependence of the errors $\bfepsilon$.
Thus the error of first kind is well controlled by all the tests in the strong case, but only by the modified versions of the Wald tests in the weak case when $N$ increases. We draw the conclusion that the modified versions are preferable to the standard ones. Table \ref{tabwaldpuis} shows that the powers of all the tests are very similar in the strong PVAR model III. The same is also true for the two modified Wald tests in the weak PVAR model IV. The empirical powers of the standard Wald tests are hardly interpretable for Model IV, because we have already seen in Table \ref{tabwald} that the standard Wald test does not well controls the error of first kind in the weak PVAR framework.

%\end{landscape}

\begin{table}[H]
 \caption{\small{Empirical size of standard and modified tests: relative frequencies (in \%)  of rejection of $H_0:\bfPhi_{22}(\nu)=0$ for $\nu=1,\dots,5$. Modified Wald Test \eqref{modwald}, which use the spectral estimator or the HAC estimator denoted: $\mathrm{{W}}_N^{\mathrm{SP}}(\nu)$ or $\mathrm{{W}}_N^{\mathrm{HAC}}(\nu)$. The number of replications is $1000$. }}
{\scriptsize
\begin{center}
\begin{tabular}{lll rrrrr rrrrr rrrrr}
\hline\hline
Model& Length $N$ & Level  & \multicolumn{5}{c}{$\mathrm{{W}}_N^{\ast}(\nu)$}& \multicolumn{5}{c}{$\mathrm{{W}}_N^{\mathrm{SP}}(\nu)$} & \multicolumn{5}{c}{$\mathrm{{W}}_N^{\mathrm{HAC}}(\nu)$}\\
& &  & \multicolumn{5}{c}{$\nu$}& \multicolumn{5}{c}{$\nu$} & \multicolumn{5}{c}{$\nu$}
\\
&&&$1$&$2$&$3$&$4$&$5$
&$\quad$ $1$&$2$&$3$&$4$&$5$& $\quad$ $1$&$2$&$3$&$4$&$5$\vspace*{0.1cm}\\
&& $\alpha=1\%$&0.8 &\textbf{0.1} &0.9& 0.8 &1.3 &$\quad$0.8& 0.7& 0.9 &1.2 &1.0&$\quad$0.8 &0.7 &1.0 &1.2 &1.1\vspace*{0.1cm}\\
I&$N=1,000$& $\alpha=5\%$& 4.5 &3.4 &4.2& 4.7 &6.2 &$\quad$5.2 &5.1 &5.3 &6.0 &\textbf{6.5}&$\quad$5.0 &5.3& 5.4 &6.1 &6.4\vspace*{0.1cm}\\
 && $\alpha=10\%$&9.6 & 8.3 &10.2 &11.5 &10.1 &$\quad$10.1 &11.1 &\textbf{12.2} &\textbf{12.1} &10.3&$\quad$9.8 &11.0 &11.9 &11.8 &10.5\vspace*{0.1cm}
 \\
 \\
 && $\alpha=1\%$&0.8& 0.4& 1.0& 1.1& 0.7 &$\quad$0.9 &1.1 &1.3& 1.0 &0.9&$\quad$1.0& 1.2 &1.3 &1.0 &0.8\vspace*{0.1cm}\\
I&$N=4,000$& $\alpha=5\%$&4.8 &4.8 &4.0& 5.6 &4.6 &$\quad$4.8 &5.6 &5.7 &6.0 &4.6&$\quad$5.0 &5.8 &5.7 &6.1 &4.5\vspace*{0.1cm}\\
 && $\alpha=10\%$&9.6 & 9.2 & 9.8 &10.5 & 9.8 &$\quad$9.8 &11.2& 11.6 &11.1 & 9.7&$\quad$10.0& 11.7 &11.4 &11.0 &10.1\vspace*{0.1cm}
 \\
\\
&& $\alpha=1\%$&\textbf{35.2}& \textbf{31.9}& \textbf{33.0} &\textbf{36.4}& \textbf{34.3}&$\quad$ \textbf{2.1}& 1.4 &\textbf{2.0}& \textbf{2.0} &1.2&$\quad$\textbf{2.1} &1.3 &1.4 &1.5 &1.1\vspace*{0.1cm}\\
II&$N=1,000$& $\alpha=5\%$&\textbf{46.6} &\textbf{44.5} &\textbf{45.5} &\textbf{51.3}& \textbf{45.6} & $\quad$6.3& 6.2 &\textbf{7.9}& \textbf{6.5}& \textbf{6.9}& $\quad$\textbf{6.6} &6.1 &\textbf{6.6} &6.1 &\textbf{6.5}\vspace*{0.1cm}\\
 && $\alpha=10\%$&\textbf{54.6} &\textbf{51.9} &\textbf{52.9} &\textbf{57.9} &\textbf{53.2}& $\quad$\textbf{12.5} &11.2& \textbf{13.2} &\textbf{12.7} &\textbf{12.2}& $\quad$11.5& 10.8 &\textbf{13.0} &\textbf{12.2}& 11.4\vspace*{0.1cm}
 \\
 \\
 && $\alpha=1\%$&\textbf{38.6}& \textbf{32.2}& \textbf{33.3}& \textbf{39.7} &\textbf{36.8}&$\quad$ 1.1& 1.4 &1.4 &1.4 &0.5& $\quad$1.0& 1.3& 1.4& 1.2 &0.5\vspace*{0.1cm}\\
II&$N=4,000$& $\alpha=5\%$&\textbf{51.9}& \textbf{46.7} &\textbf{46.3}& \textbf{52.0} &\textbf{48.6}  & $\quad$5.9 &5.6& 5.2& 5.8 &4.4& $\quad$6.0& 5.5 &4.8 &6.0 &4.2\vspace*{0.1cm}\\
 && $\alpha=10\%$&\textbf{57.5} &\textbf{54.1}& \textbf{53.1}& \textbf{57.5}& \textbf{55.1}& $\quad$11.1& 11.1& 10.5 &11.9 & 9.6& $\quad$10.4 &10.4 &10.8 &11.8 & 9.0\vspace*{0.1cm}
 \\
\hline\hline \multicolumn{17}{l}{I: Strong PVAR$(1)$ model
(\ref{dgp})-(\ref{bruitfort}) with unknown parameter given in Table \ref{param1}.} \\
\multicolumn{17}{l}{II: Weak PVAR$(1)$ model
(\ref{dgp})-(\ref{bruitfaible}) with unknown parameter given in Table \ref{param1}.}\\
\end{tabular}
\end{center}
}
\label{tabwald}
\end{table}

\begin{table}[H]
\caption{Parameters of DGP models used in the simulation to test $H_0:\bfPhi_{22}(\nu)=0$ for $\nu=1,\dots,5$.}\label{param1}
\vspace{0.1cm}
\centering
\begin{tabular}{l cc cccc cc cccc cc cccc cc cccc cc cccc}
  \toprule
  MODEL & & \multicolumn{2}{c}{$\bfPhi(1)$} & & \multicolumn{2}{c}{$\bfPhi(2)$}& & \multicolumn{2}{c}{$\bfPhi(3)$}& & \multicolumn{2}{c}{$\bfPhi(4)$}& & \multicolumn{2}{c}{$\bfPhi(5)$}\\
  \midrule
  DGP & &-1.43& 0.00&& 0.46& 0.00&&1.23& 0.00&&0.30& 0.00&&0.90& 0.00\\
          & &0.00& 0.00&& 0.00& 0.00&&0.00& 0.00&&0.00& 0.00&&0.00& 0.00\\
  \bottomrule
\end{tabular}
\end{table}

\begin{table}[H]
 \caption{\small{Empirical power of standard and modified tests: relative frequencies (in \%)  of rejection of $H_0:\bfPhi_{22}(\nu)=0$ for $\nu=1,\dots,5$. Modified Wald Test \eqref{modwald}, which use the spectral estimator or the HAC estimator denoted: $\mathrm{{W}}_N^{\mathrm{SP}}(\nu)$ or $\mathrm{{W}}_N^{\mathrm{HAC}}(\nu)$. The number of replications is $1000$. }}
{\scriptsize
\begin{center}
\begin{tabular}{lll rrrrr rrrrr rrrrr}
\hline\hline
Model& Length $N$ & Level  & \multicolumn{5}{c}{$\mathrm{{W}}_N^{\ast}(\nu)$}& \multicolumn{5}{c}{$\mathrm{{W}}_N^{\mathrm{SP}}(\nu)$} & \multicolumn{5}{c}{$\mathrm{{W}}_N^{\mathrm{HAC}}(\nu)$}\\
& &  & \multicolumn{5}{c}{$\nu$}& \multicolumn{5}{c}{$\nu$} & \multicolumn{5}{c}{$\nu$}
\\
&&&$1$&$2$&$3$&$4$&$5$
&$\quad$ $1$&$2$&$3$&$4$&$5$& $\quad$ $1$&$2$&$3$&$4$&$5$\vspace*{0.1cm}\\
&& $\alpha=1\%$&79.0& 100.0 & 44.8 & 51.4 & 54.7 &$\quad$ 78.8 &100.0  &49.2  &52.0 & 55.2&$\quad$78.8 &100.0  &48.7 & 52.1 & 54.7\vspace*{0.1cm}\\
III&$N=4,000$& $\alpha=5\%$& 93.3 &100.0 & 67.9 & 75.1 & 76.4 &$\quad$93.1 &100.0 & 70.1 & 75.1 & 76.9  &$\quad$93.0& 100.0  &70.3  &75.2 & 76.6\vspace*{0.1cm}\\
 && $\alpha=10\%$&97.1& 100.0 & 77.2 & 82.7 & 84.8 &$\quad$97.4 &100.0  &78.7  &83.0  &85.0&$\quad$97.3 &100.0 & 78.5 & 83.2 & 84.9\vspace*{0.1cm}
 \\
\\
&& $\alpha=1\%$&59.3 &85.4 &52.3 &54.5 &54.0&$\quad$  7.8 &30.7  &6.5  &5.2 & 5.2&$\quad$7.4& 30.3 & 5.8 & 5.4 & 5.1\vspace*{0.1cm}\\
IV&$N=4,000$& $\alpha=5\%$&68.4 &90.8 &63.6& 66.3 &64.1& $\quad$23.4 &53.4 &17.7 &16.9 &15.4& $\quad$22.4 &53.1 &17.7 &16.2 &15.0\vspace*{0.1cm}\\
 && $\alpha=10\%$&74.4 &92.9 &69.3 &71.7& 71.7& $\quad$32.8 &65.6 &27.1 &25.7& 24.3& $\quad$32.0 &66.0 &27.4 &25.2 &24.2\vspace*{0.1cm}
 \\
\hline\hline \multicolumn{17}{l}{III: Strong PVAR$(1)$ model
(\ref{dgp})-(\ref{bruitfort}) with unknown parameter given in Table \ref{param2}.} \\
\multicolumn{17}{l}{IV: Weak PVAR$(1)$ model
(\ref{dgp})-(\ref{bruitfaible}) with unknown parameter given in Table \ref{param2}.}\\
\end{tabular}
\end{center}
}
\label{tabwaldpuis}
\end{table}

\begin{table}[H]
\caption{Parameters of DGP models used in the simulation of empirical power to test $H_0:\bfPhi_{22}(\nu)=0$ for $\nu=1,\dots,5$.}\label{param2}
\vspace{0.1cm}
\centering
\begin{tabular}{l cc cccc cc cccc cc cccc cc cccc cc cccc}
  \toprule
  MODEL & & \multicolumn{2}{c}{$\bfPhi(1)$} & & \multicolumn{2}{c}{$\bfPhi(2)$}& & \multicolumn{2}{c}{$\bfPhi(3)$}& & \multicolumn{2}{c}{$\bfPhi(4)$}& & \multicolumn{2}{c}{$\bfPhi(5)$}\\
  \midrule
  DGP & &-1.43& 0.00&& 0.46& 0.00&&1.23& 0.00&&0.30& 0.00&&0.90& 0.00\\
          & &0.00& 0.05&& 0.00& 0.05&&0.00& 0.05&&0.00& 0.05&&0.00& 0.05\\
  \bottomrule
\end{tabular}
\end{table}

%%% ----------------------------------------------------------------------
%%% ----------------------------------------------------------------------
\section{Application to real data}\label{real}
\noindent In this section, we consider the daily returns of two European stock market indices: CAC 40 (Paris) and DAX (Frankfurt), from March $3$, $1990$ to March $10$, $2022$. The data were obtained from \textit{Yahoo Finance}. Because of the legal holidays, many weeks comprise less than five observations. We preferred removing the entire weeks when there was less than five data available, giving a bivariate time series of sample size equal to $7060$. The period $\nu=5$ is naturally selected.

\begin{figure}[H]
  \centering
  \includegraphics[width=1\linewidth]{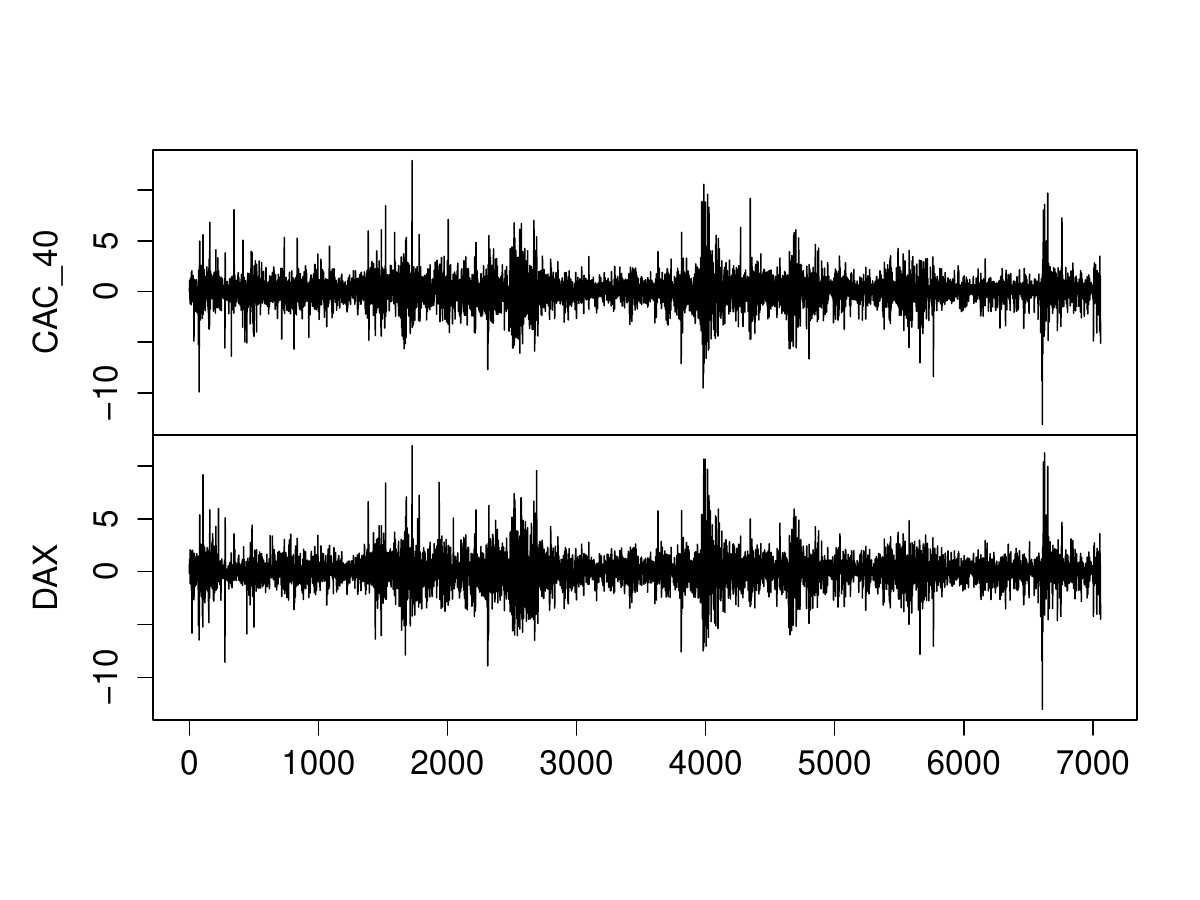}
  \caption{Log-returns of CAC 40 (Paris) and DAX (Frankfurt)}
  \label{fig:log_returns}
\end{figure}

\noindent In order to analyse these two European indices, we fitted a PVAR model of order $1$ to the bivariate series of observations:
$$\bfY_{ns+\nu} = \bfPhi(\nu)\bfY_{ns+\nu-1} +\bfepsilon_{ns+\nu}\quad \nu = 1,\ldots,5,$$
where
$\bfY_t = \left(r_t^1,r_t^2\right)^\top$ and $r_t^1$, $r_t^2$ represents the log-return of CAC 40 and DAX respectively.
The log-return is defined as $r_t = 100\times\ln{(I_t/I_{t-1})}$ where $I_t$ represents the value of the index at time $t$. Seasonal means are first removed from the series, meaning that a model is formulated by examining $\bfY_{ns+\nu}-\bfmu(\nu)$. The two time series of log-returns are displayed in Figure~\ref{fig:log_returns}.

\noindent We present in Table~\ref{table:results} the estimated parameters $\hat{\bfbeta}=(\hat{\bfbeta}(1),\dots,\hat{\bfbeta}(5))^\top$ and their estimated standard error proposed in the strong case~(see Remark \ref{rem1}) denoted $\hat{\sigma}_\mathrm{S}$ and the weakly consistent estimators proposed~(\ref{estThetaSP}) and~(\ref{estThetaHAC}), denoted respectively  by $\hat{\sigma}_\mathrm{{SP}}$ and $\hat{\sigma}_\mathrm{{HAC}}$; $\hat{\bfSigma}(\nu)$ represents the estimated variance of residuals $\hat{\bfepsilon}(\nu)$. The $p$-values of the $t$-statistic of  $\hat{\bfbeta}$ and those of the standard and modified versions of the Wald tests are denoted: pval$_S$, pval$_\mathrm{{SP}}$,  pval$_\mathrm{{HAC}}$, pval$^\mathrm{{W}}_S$, pval$^\mathrm{{W}}_\mathrm{{SP}}$ and  pval$^\mathrm{{W}}_\mathrm{{HAC}}$, where the exponent W stands for Wald.
The $p$-values less than 5\% are in bold, those less than 1\% are underlined. The autoregressive coefficients $\hat{\bfPhi}_{ij}(\nu)$ for $i,j=1,2$ are rather small on Monday, Tuesday and Friday. Four of them are significant at the 1\% level in the strong case on Wednesday and Thursday. In the weak case, none of them are significant at the 1\% level. This is in accordance with the results of \cite{frs11} who showed that the log-returns of these two European stock market indices constitute the weak periodic white noises. As in~\cite{frs11}, the estimated variance of residuals shows that, the volatility is considerably greater on Monday and smaller for the other days.

\begin{table}[H]
\small
\begin{center}
\caption{\label{table:results}Least squares estimators used to fit the log-returns of CAC 40 and DAX data to a bivariate PVAR model with $\nu=5$; the $\hat{\sigma}_\mathrm{S}$, $\hat{\sigma}_\mathrm{{SP}}$ and $\hat{\sigma}_\mathrm{{HAC}}$ represent the standard errors in the strong case and for our proposed estimators in the weak case; pval$_S$, pval$_\mathrm{{SP}}$,  pval$_\mathrm{{HAC}}$, pval$^\mathrm{{W}}_S$, pval$^\mathrm{{W}}_\mathrm{{SP}}$ and  pval$^\mathrm{{W}}_\mathrm{{HAC}}$ correspond to the $p$-values of the $t$-statistic of  $\hat{\bfbeta}$ and those of the standard and modified versions of the Wald tests are also presented; $\hat{\bfSigma}_{\bfepsilon}(\nu)$ represents the estimated variance of residuals $\hat{\bfepsilon}(\nu)$. The $p$-values less than 5\% are in bold, those less than 1\% are underlined.}
\centering
\begin{tabular}{ccccccccc ccc}
  \hline
 & $\hat{\bfbeta}$ & $\hat{\sigma}_\mathrm{S}$ & $\hat{\sigma}_\mathrm{{SP}}$ & $\hat{\sigma}_\mathrm{{HAC}}$ & pval$_S$& pval$_\mathrm{{SP}}$ & pval$_\mathrm{{HAC}}$& pval$^\mathrm{{W}}_S$& pval$^\mathrm{{W}}_\mathrm{{SP}}$ & pval$^\mathrm{{W}}_\mathrm{{HAC}}$  & $\vec{\left(\hat{\bfSigma}_{\bfepsilon}(\nu)\right)}$\\
  \hline
  1 & -0.0349 & 0.0707 & 0.1456 & 0.1078 & 0.6220 & 0.8107 & 0.7464 & 0.6219 & 0.8107 & 0.7464 & 3.5457 \\
  2 & 0.0153 & 0.0731 & 0.1480 & 0.1152 & 0.8346 & 0.9180 & 0.8947 & 0.8346 & 0.9179 & 0.8947 & 3.2404 \\
  3 & -0.0070 & 0.0706 & 0.0992 & 0.0862 & 0.9215 & 0.9441 & 0.9357 & 0.9215 & 0.9441 & 0.9357 & 3.2404 \\
  4 & -0.0378 & 0.0729 & 0.1019 & 0.1004 & 0.6044 & 0.7106 & 0.7066 & 0.6043 & 0.7106 & 0.7065 & 3.7859 \\
  5 & -0.0506 & 0.0399 & 0.0524 & 0.0486 & 0.2045 & 0.3339 & 0.2975 & 0.2043 & 0.3337 & 0.2973 & 1.7366 \\
  6 & -0.0270 & 0.0420 & 0.0663 & 0.0622 & 0.5214 & 0.6843 & 0.6651 & 0.5213 & 0.6842 & 0.6650 & 1.5938 \\
  7 & -0.0020 & 0.0386 & 0.0551 & 0.0474 & 0.9591 & 0.9713 & 0.9667 & 0.9591 & 0.9713 & 0.9667 & 1.5938 \\
  8 & -0.0246 & 0.0407 & 0.0742 & 0.0624 & 0.5450 & 0.7399 & 0.6932 & 0.5449 & 0.7398 & 0.6931 & 1.9297 \\
  9 & -0.3001 & 0.0532 & 0.1296 & 0.1200 & \underline{\textbf{0.0000}} & \textbf{0.0207} & \textbf{0.0125} & \underline{\textbf{0.0000}} & \textbf{0.0206} & \textbf{0.0124} & 1.6817 \\
  10 & -0.1256 & 0.0545 & 0.0781 & 0.0693 & \textbf{0.0214} & 0.1078 & 0.0704 & \textbf{0.0213} & 0.1076 & 0.0702 & 1.4736 \\
  11 & 0.2605 & 0.0505 & 0.1381 & 0.1265 & \underline{\textbf{0.0000}} & 0.0595 & \textbf{0.0396} & \underline{\textbf{0.0000}} & 0.0592 & \textbf{0.0394} & 1.4736 \\
  12 & 0.0360 & 0.0517 & 0.0607 & 0.0642 & 0.4869 & 0.5535 & 0.5752 & 0.4868 & 0.5534 & 0.5751 & 1.7671 \\
  13 & -0.1498 & 0.0548 & 0.1020 & 0.0844 & \underline{\textbf{0.0063}} & 0.1422 & 0.0761 & \underline{\textbf{0.0062}} & 0.1420 & 0.0758 & 2.0261 \\
  14 & -0.0744 & 0.0551 & 0.0639 & 0.0715 & 0.1767 & 0.2445 & 0.2979 & 0.1765 & 0.2443 & 0.2977 & 1.7907 \\
  15 & 0.1862 & 0.0538 & 0.1021 & 0.0855 & \underline{\textbf{0.0006}} & 0.0683 & \textbf{0.0295} & \underline{\textbf{0.0005}} & 0.0681 & \textbf{0.0294} & 1.7907 \\
  16 & 0.0955 & 0.0541 & 0.0620 & 0.0715 & 0.0778 & 0.1240 & 0.1824 & 0.0776 & 0.1237 & 0.1822 & 2.0471 \\
  17 & -0.0227 & 0.0521 & 0.0681 & 0.0647 & 0.6627 & 0.7385 & 0.7252 & 0.6626 & 0.7384 & 0.7251 & 1.7824 \\
  18 & -0.0055 & 0.0522 & 0.0670 & 0.0654 & 0.9156 & 0.9343 & 0.9326 & 0.9156 & 0.9343 & 0.9326 & 1.5118 \\
  19 & 0.0694 & 0.0520 & 0.0739 & 0.0699 & 0.1823 & 0.3480 & 0.3209 & 0.1821 & 0.3479 & 0.3207 & 1.5118 \\
  20 & 0.0420 & 0.0521 & 0.0734 & 0.0705 & 0.4202 & 0.5677 & 0.5517 & 0.4201 & 0.5676 & 0.5516 & 1.7861 \\
   \hline
\end{tabular}
\end{center}
\end{table}
%%%------------------------------------------------------------------------------------------------
%%%------------------------------------------------------------------------------------------------
\section{Conclusions}\label{conclusion}
\noindent In this work, we have established under mild assumptions, the asymptotic distribution of the least squares estimator of the model parameters in PVAR time series models with dependent but uncorrelated errors. Our results extend Theorem~1 of~\cite{UD09} for PVAR models with independent errors. Note that if $s=1$, we retrieve the result on weak VAR obtained by~\cite{FR07}. The asymptotic covariance matrix of the least squares estimators obtained under independent errors is generally not consistent in the weak PVAR case. For statistical inference problem, including in particular, the significance tests on the parameters, the assumption of independent errors can be quite misleading when analysing data from PVAR models with dependent errors.

We proposed two estimators of the asymptotic variance matrix: the spectral density estimator and the heteroskedasticity and autocorrelation consistent estimator based on kernel methods. The empirical results of Sections~\ref{simul} and~\ref{real} illustrate the applicability of our theoretical results using a consistent estimator of the asymptotic variance matrix of the least square estimators of weak PVAR parameters. In future works, we intend to study how the existing identification and diagnostic checking (see e.g.~\cite{UD09}) procedures should be adapted in the weak PVAR framework considered in the present paper. The asymptotic covariance matrix of the least squares estimators of a weak PVAR model is no longer block diagonal with respect to seasons and depends on the fourth-order moments of the innovation process (through
the matrix $\bfPsi(\nu)$).

%% The Appendices part is started with the command \appendix;
%% appendix sections are then done as normal sections
\appendix
\section{Appendix : Proofs of the main results}\label{app}
The proof of Theorem \ref{thm} is quite technical. This is adaptation of the arguments used in \cite{frs11}.
\subsection{Proof of Theorem \ref{thm}}
%\begin{proof}

The proof is quite long so we divide it in several steps.
\paragraph{$\diamond$ Step 1: preliminaries} \ \\

In view of \eqref{MAinf}, it is easy to see that $\bfX_n^{\top}(\nu)$
is a measurable function of the random vectors
$\{ \bfepsilon_{ns+\nu-k}, k \geq 1 \}$. Thus the assumption {\bf (A2)} of the error term
$(  \bfepsilon_n^{\ast} )_{n\in\mathbb{Z}}$ allows us to show that
$( \vec\{\bfepsilon_{ns+\nu}\bfX_n^{\top}(\nu) \} )_{n\in\mathbb{Z}}$
is a stationary and ergodic  sequence. Applying the ergodic theorem, we obtain that
\begin{equation}\label{ortho}
  N^{-1}\sum_{n=0}^{N-1} \vec\{ \bfepsilon_{ns+\nu}\bfX_n^{\top}(\nu) \} \cas\mathbb{E}\left[\vec\{ \bfepsilon_{ns+\nu}\bfX_n^{\top}(\nu) \}\right]= \bfzero,
\end{equation}
by using the non-correlation between $\bfepsilon_{ns+\nu}$'s (see {\bf (A0)}) and where
$\bfzero$ is the $\{ d^2 p(\nu) \} \times 1$ null vector.
\paragraph{$\diamond$ Step 2: convergence in distribution of $N^{-1/2}\sum_{n=0}^{N-1} \vec\{ \bfepsilon_{ns+\nu}\bfX_n^{\top}(\nu) \}$} \ \\

Using the stationarity of $( \vec\{\bfepsilon_{ns+\nu}\bfX_n^{\top}(\nu) \} )_{n\in\mathbb{Z}}$, we have
\begin{align*}
\var&\left\lbrace \frac{1}{\sqrt{N}}\sum_{n=0}^{N-1} \vec\{ \bfepsilon_{ns+\nu}\bfX_n^{\top}(\nu) \}\right\rbrace \\&=\frac{1}{N}\sum_{n=0}^{N-1}\sum_{n'=0}^{N-1}\cov\left\lbrace \vec\{ \bfepsilon_{ns+\nu}\bfX_n^{\top}(\nu) \},\vec\{ \bfepsilon_{n's+\nu}\bfX_{n'}^{\top}(\nu) \}\right\rbrace \\
&=\frac{1}{N}\sum_{h=-N+1}^{N-1}\left( N-|h|\right)c(\nu,h),
\end{align*}
where \[c(\nu,h)=\cov\left( \vec\{ \bfepsilon_{ns+\nu}\bfX_n^{\top}(\nu) \},\vec\{ \bfepsilon_{(n-h)s+\nu}\bfX_{n-h}^{\top}(\nu) \}\right).\]
By the dominated convergence theorem, it follows that
$$\bfPsi(\nu)=\sum_{h=-\infty}^{\infty}\cov\left( \vec\{ \bfepsilon_{ns+\nu}\bfX_n^{\top}(\nu) \},\vec\{ \bfepsilon_{(n-h)s+\nu}\bfX_{n-h}^{\top}(\nu) \}\right). $$
The existence of the last sum is a consequence of {\bf (A3)} and the \cite{D68} inequality. Using \eqref{ortho} and the elementary relations $\vec(ab^\top)=b\otimes a$ for any vectors $a$ and $b$, and $(A\otimes B)(C\otimes D)= (AC)\otimes (BD)$ for matrices of appropriate sizes (see \cite{Lu05}), it follows that
$$  \bfPsi(\nu)= \sum_{h=-\infty}^{\infty}\mathbb{E}\left(\bfX_n(\nu)\bfX_{n-h}^{\top}(\nu) \otimes \bfepsilon_{ns+\nu}\bfepsilon_{(n-h)s+\nu}^{\top}\right).$$

Let $\bfepsilon_n(\nu) = (\bfepsilon_{ns+\nu-1}^{\top},\ldots,\bfepsilon_{ns+\nu-p(\nu)}^{\top})^{\top}$, $n=0,1,\ldots,N-1$,
be a $\{ d p(\nu) \} \times 1$ random vectors. In the sequel, we need the elementary identity $\vec(ABC)=(I\otimes AB)\vec(C)$ (see \cite{Lu05}). In view of \eqref{MAinf}, we have for all $r\geq0$
\begin{align}\nonumber
\vec\{ \bfepsilon_{ns+\nu}\bfX_n^{\top}(\nu) \}&=
%\sum_{i=0}^{\infty} \vec\{ \bfepsilon_{ns+\nu}\bfepsilon_{n-i}^{\top}(\nu)\bfC_i^{\top}(\nu) \}=
\sum_{i=0}^{\infty}\left(\bfI_{dp(\nu)}\otimes\bfepsilon_{ns+\nu}\bfepsilon_{n-i}^{\top}(\nu)\right)\vec\left(\bfI_{p(\nu)}\otimes\bfC_i^{\top}(\nu) \right)
\\ \label{sumfini}&=\bfW_{n,r}(\nu)+\bfU_{n,r}(\nu),
\end{align}
where
\begin{align*}
\bfW_{n,r}(\nu)&=\sum_{i=0}^{r}\left(\bfI_{dp(\nu)}\otimes\bfepsilon_{ns+\nu}\bfepsilon_{n-i}^{\top}(\nu)\right)\vec\left(\bfI_{p(\nu)}\otimes\bfC_i^{\top}(\nu) \right)
\\\bfU_{n,r}(\nu)&=\sum_{i=r+1}^{\infty}\left(\bfI_{dp(\nu)}\otimes\bfepsilon_{ns+\nu}\bfepsilon_{n-i}^{\top}(\nu)\right)\vec\left(\bfI_{p(\nu)}\otimes\bfC_i^{\top}(\nu) \right).
\end{align*}
The processes $(\bfW_{n,r}(\nu))_{n\in\mathbb{Z}}$ and $(\bfU_{n,r}(\nu))_{n\in\mathbb{Z}}$ are stationary and centered.
Moreover, under Assumption {\bf (A3)} and  $r$ fixed, the process $(\bfW_{n,r}(\nu))_{n\in\mathbb{Z}}$ is strongly mixing (see Theorem 14.1 in \cite{D94}), with mixing coefficients
$\alpha_{\bfW_r}(h)\leq \alpha_{\bfepsilon}\left(\max\{0,h-1\}\right)$. Thus {\bf (A3)} implies $\sum_{h=0}^{\infty}\{\alpha_{\bfW_r}(h)\}^{\kappa/(2+\kappa)}< \infty$ and using the H\"{o}der inequality, we obtain that $\|\bfW_{n,r}(\nu)\|_{2+\kappa}< \infty$ for some $\kappa>0$.  The central limit theorem for strongly mixing processes (see \cite{herr})  implies that  $N^{-1/2}\sum_{n=0}^{N-1}\bfW_{n,r}(\nu)$ has a limiting $\mathcal{N}(0,\bfPsi_r(\nu))$ distribution with
$$\bfPsi_r(\nu)=\lim_{N\rightarrow\infty}\var\left( \frac{1}{\sqrt{N}}\sum_{n=0}^{N-1}\bfW_{n,r}(\nu)\right)=\sum_{h=-\infty}^{\infty}
\cov\left(\bfW_{n,r}(\nu),\bfW_{n-h,r}(\nu)\right) .$$
Since $ N^{-1/2}\sum_{n=0}^{N-1}\bfW_{n,r}(\nu)$ and $ N^{-1/2}\sum_{n=0}^{N-1}\vec\{ \bfepsilon_{ns+\nu}\bfX_n^{\top}(\nu)\} $ have zero expectation, we shall have
$$\lim_{r\rightarrow\infty}\var\left( \frac{1}{\sqrt{N}}\sum_{n=0}^{N-1}\bfW_{n,r}(\nu)\right)=\var\left\lbrace \frac{1}{\sqrt{N}}\sum_{n=0}^{N-1} \vec\{ \bfepsilon_{ns+\nu}\bfX_n^{\top}(\nu) \}\right\rbrace  ,$$
 as soon as
 \begin{equation}
\label{Zt}
\lim_{r\to\infty}\limsup_{N\to\infty}\mathbb{P}\left\{\left\|\frac{1}{\sqrt{N}}\sum_{n=0}^{N-1}\bfU_{n,r}(\nu)\right\|>\varepsilon\right\}=0
\end{equation}
for every $\varepsilon>0$.
%\begin{equation}\label{conv_L2}
%\lim_{r\rightarrow\infty}\mathbb{E}\left[ \left\| \frac{1}{\sqrt{N}}\sum_{n=0}^{N-1} \vec\{ \bfepsilon_{ns+\nu}\bfX_n^{\top}(\nu) \}- \frac{1}{\sqrt{N}}\sum_{n=0}^{N-1}\bfW_{n,r}(\nu)\right\|^2\right] =0 .
%\end{equation}
As a consequence we will have $\lim_{r\rightarrow\infty}\bfPsi_r(\nu)=\bfPsi(\nu)$. The result \eqref{Zt} follows from a
straightforward adaptation of Theorem 7.7.1 and Corollary 7.7.1 of Anderson (see \cite{anderson} pages 425-426). Indeed, by stationarity we have
\begin{eqnarray*}
\var\left(\frac{1}{\sqrt{N}}\sum_{n=0}^{N-1}\bfU_{n,r}(\nu)\right)&=&\frac{1}{N}\sum_{n,n'=0}^{N-1}\cov\left(\bfU_{n,r}(\nu),\bfU_{n',r}(\nu)\right)\\
&=&\frac{1}{N}\sum_{|h|<N-1}(N-|h|)\cov\left(\bfU_{n,r}(\nu),\bfU_{n-h,r}(\nu)\right)\\
&\leq&\sum_{h=-\infty}^{\infty}\left\|\cov\left(\bfU_{n,r}(\nu),\bfU_{n-h,r}(\nu)\right)\right\|.
\end{eqnarray*}
Because $\|\bfC_{i}\|\leq K\rho^{i}$ for $\rho\in[0,1[$ and $K>0$ and in view of \eqref{sumfini}, we have
$$\left\|\bfU_{n,r}(\nu)\right\|\leq K\sum_{i=r+1}^{\infty}\rho^{i}\|\bfepsilon_{ns+\nu}\|\|\bfepsilon_{n-i}(\nu)\|.$$
Under {\bf (A3)} we have $\mathbb{E}\|\bfepsilon_{ns+\nu}\|^{4+2\kappa}<\infty$, it follows from the Hölder inequality that
\begin{equation}
\label{firstmajo}
\sup_h\left\|\cov\left(\bfU_{n,r}(\nu),\bfU_{n-h,r}(\nu)\right)\right\|=\sup_h\left\|\mathbb{E}\left(\bfU_{n,r}(\nu)\bfU_{n-h,r}^\top(\nu)\right)\right\|\leq K\rho^r.
\end{equation}
Let $h>0$ such that $[h/2]>r$. Write
$$\bfU_{n,r}(\nu)=\bfU_{n,r}^{h^-}(\nu)+\bfU_{n,r}^{h^+}(\nu),$$
where
\begin{align*}
\bfU_{n,r}^{h^-}(\nu)&=\sum_{i=r+1}^{[h/2]}\left(\bfI_{dp(\nu)}\otimes\bfepsilon_{ns+\nu}\bfepsilon_{n-i}^{\top}(\nu)\right)\vec\left(\bfI_{p(\nu)}\otimes\bfC_i^{\top}(\nu) \right),
\\
 \bfU_{n,r}^{h^+}(\nu)&=\sum_{i=[h/2]+1}^{\infty}\left(\bfI_{dp(\nu)}\otimes\bfepsilon_{ns+\nu}\bfepsilon_{n-i}^{\top}(\nu)\right)\vec\left(\bfI_{p(\nu)}\otimes\bfC_i^{\top}(\nu) \right).
\end{align*}
Note that $\bfU_{n,r}^{h^-}(\nu)$ belongs to the $\sigma$-field generated by $\{\bfepsilon_{ns+\nu},\bfepsilon_{ns+\nu-1},\dots,\bfepsilon_{ns+\nu-[h/2]}\}$ and
that $\bfU_{n-h,r}(\nu)$ belongs to the $\sigma$-field generated by $\{\bfepsilon_{(n-h)s+\nu},\bfepsilon_{(n-h-1)s+\nu-1},,\dots\}$.
By {\bf (A3)}, $\mathbb{E}\|\bfU_{n,r}^{h^-}(\nu)\|^{2+\kappa}<\infty$
and $\mathbb{E}\|\bfU_{n-h,r}(\nu)\|^{2+\kappa}<\infty$. Davydov's inequality (see \cite{D68}) then entails that
\begin{equation}
\label{secondmajo}
\left\|\cov\left(\bfU_{n,r}^{h^-}(\nu),\bfU_{n-h,r}(\nu)\right)\right\|\leq K\alpha_{\bfepsilon}^{\kappa/(2+\kappa)}([h/2]).
\end{equation}
By the argument used to show (\ref{firstmajo}), we also have
\begin{equation}
\label{thirdmajo}
\left\|\cov\left(\bfU_{n,r}^{h^+}(\nu),\bfU_{n-h,r}(\nu)\right)\right\|\leq K\rho^h\rho^r.
\end{equation}
In view of (\ref{firstmajo}), (\ref{secondmajo}) and (\ref{thirdmajo}), we have
$$\sum_{h=0}^{\infty}\left\|\cov\left(\bfU_{n,r}(\nu),\bfU_{n-h,r}(\nu)\right)\right\|\leq
K\rho^r+K\sum_{h=r}^{\infty}\alpha_{\bfepsilon}^{\kappa/(2+\kappa)}(h)\to 0$$
as $r\to\infty$ by {\bf (A3)}. We have the same bound for $h<0$.
This implies that
\begin{align}\label{koko}
& \sup_N\var\left(\frac{1}{\sqrt{N}}\sum_{n=0}^{N-1}\bfU_{n,r}(\nu)\right)\xrightarrow[r\to\infty]{} 0.
\end{align}
The conclusion of \eqref{Zt} follows from the Markov inequality.
% which completes the proof of \eqref{Zt}.

From a standard result (see {\em e.g.} Proposition 6.3.9 in \cite{BD91}), we deduce that
$$\frac{1}{\sqrt{N}}\sum_{n=0}^{N-1}\vec\{ \bfepsilon_{ns+\nu}\bfX_n^{\top}(\nu) \}=
\frac{1}{\sqrt{N}}\sum_{n=0}^{N-1}\bfW_{n,r}(\nu)+\frac{1}{\sqrt{N}}\sum_{n=0}^{N-1}\bfU_{n,r}(\nu)\stackrel{d}{\to}{\cal N}\left(0,\bfPsi(\nu)\right),$$
which completes the proof of \eqref{th1a}.
\paragraph{$\diamond$ Step 3: existence and invertibility of the matrix $\bfOmega(\nu)$} \ \\

By ergodicity of the centred  process $(\bfX_n(\nu))_{n\in\mathbb{Z}}\in\mathbb{R}^{ d p(\nu) }$, we deduce that
\begin{equation}\label{omega}
\frac{1}{N} \bfX_n(\nu) \bfX_n^{\top}(\nu) \cas\bfOmega(\nu):=\mathbb{E}\left(\bfX_n(\nu)\bfX_n^{\top}(\nu)\right).
\end{equation}
From \eqref{sumfini} we obtain that
\begin{align*}
\mathbb{E}\left(\bfX_n(\nu)\bfX_n^{\top}(\nu)\right)&=\mathbb{E}\left[ \left( \sum_{i=0}^{\infty}\left(\bfI_{p(\nu)}\otimes\bfC_i(\nu) \right)\bfepsilon_{n-i}(\nu)\right) \left(\sum_{j=0}^{\infty}\left(\bfI_{p(\nu)}\otimes\bfC_j(\nu) \right)\bfepsilon_{n-j}(\nu)\right)^\top \right] \\
&=\sum_{i=0}^{\infty}\sum_{j=0}^{\infty} \left(\bfI_{p(\nu)}\otimes\bfC_i(\nu) \right)\mathbb{E}\left[\bfepsilon_{n-i}(\nu)\bfepsilon_{n-j}^\top(\nu) \right]\left(\bfI_{p(\nu)}\otimes\bfC_j^\top(\nu) \right)
 \\
&=\sum_{i=0}^{\infty} \left(\bfI_{p(\nu)}\otimes\bfC_i(\nu) \right)
\left(\bfI_{p(\nu)}\otimes\bfSigma_{\bfepsilon}(\nu) \right)\left(\bfI_{p(\nu)}\otimes\bfC_i^\top(\nu) \right)
\\
& \leq K\, \sum_{i\geq 0}\rho^i
<\infty .
\end{align*}
Therefore the matrix $\bfOmega(\nu)$ exists almost surely.

If the matrix $\bfOmega(\nu)$ is not invertible, there exists some real constants $c_1,\dots,c_{dp(\nu)}$ not all equal to zero such that $\mathbf{c}^\top\bfOmega(\nu)\mathbf{c}=0$, where $\mathbf{c}=(c_1,\dots,c_{dp(\nu)})^\top$.
For $i=1,\dots,dp(\nu)$, let $\bfX_{i,n}(\nu)$ be the $i$-th component of $\bfX_{n}(\nu)$ and denotes by $\bfOmega_{ji}(\nu)$
 the $(i,j)$-th component of $\bfOmega(\nu)$.
We obtain that
\begin{align*}
\sum_{i=1}^{dp(\nu)}\sum_{j=1}^{dp(\nu)}c_j\bfOmega_{ji}(\nu)c_i=\sum_{i=1}^{dp(\nu)}\sum_{j=1}^{dp(\nu)}\mathbb{E}\left[ \left( c_j\bfX_{j,n}(\nu)\right) \left(  c_i\bfX_{i,n}(\nu)\right) \right]=
\mathbb{E}\left[ \left( \sum_{k=1}^{dp(\nu)}c_k\bfX_{k,n}(\nu)\right)^2 \right]&=0,
\end{align*}
which implies that
\begin{align*}%\label{hypInvJ}
\sum_{k=1}^{dp(\nu)}c_k\bfX_{k,n}(\nu)&=0 \ \ \mathrm{a.s.}
\text{ or equivalenty }\ \ \mathbf{c}^\top\bfX_n(\nu)=\sum_{i=0}^{\infty}\mathbf{c}^\top\left(\bfI_{p(\nu)}\otimes\bfC_i(\nu) \right)\bfepsilon_{n-i}(\nu)=0 \ \ \mathrm{a.s.}
\end{align*}
This is in contradiction with the assumption that $\bfSigma_{\bfepsilon}(\nu)$ is not equal to zero. Therefore $\mathbf{c}^\top\bfX_n(\nu)$ is not almost surely equal to zero and $\bfOmega(\nu)$ is almost surely invertible.
\paragraph{$\diamond$ Step 4: convergence in probability of $\hat{\bfbeta}(\nu)$} \ \\

Using the relation~(\ref{hatBnu}),
we can write:
\[
   \hat{\bfB}(\nu) - \bfB(\nu) = N^{-1} \bfE(\nu) \bfX^{\top}(\nu) \{ N^{-1} \bfX(\nu) \bfX^{\top}(\nu) \}^{-1}.
\]
Noting that
$\sum_{n=0}^{N-1} \vec\{ \bfepsilon_{ns+\nu} \bfX_n^{\top}(\nu) \} = \vec\{ \bfE(\nu)\bfX^{\top}(\nu) \}$, from \eqref{th1a},
it follows that
$N^{-1/2}\vec\{ \bfE(\nu) \bfX^{\top}(\nu) \}
 \cd N_{d^2p(\nu)}(\bfzero, \bfPsi(\nu) )$. Applying the ergodic theorem and from \eqref{ortho}, we have
$N^{-1}\vec\{ \bfE(\nu) \bfX^{\top}(\nu) \} \cas \bfzero$,
where the dimension of
$\bfzero$
is
$\{ d^2 p(\nu) \} \times 1$,
and also
$\{ N^{-1}\bfX(\nu)\bfX^{\top}(\nu) \}^{-1} \cas \bfOmega^{-1}(\nu)$;
these results show~(\ref{th1b}).
\paragraph{$\diamond$ Step 5: convergence in distribution of $N^{1/2}\{ \hat{\bfbeta}(\nu) - \bfbeta(\nu) \}$} \ \\

Since
\begin{align}
\nonumber
  N^{1/2}\{ \hat{\bfbeta}(\nu) - \bfbeta(\nu) \} &=
  \left[ \{ N^{-1}\bfX(\nu) \bfX^{\top}(\nu) \}^{-1} \otimes \bfI_d \right] N^{-1/2} \{ \bfX(\nu) \otimes \bfI_d \} \bfe(\nu),\\ \label{betahatnu}
 & =
  \left[ \{ N^{-1}\bfX(\nu) \bfX^{\top}(\nu) \}^{-1} \otimes \bfI_d \right] N^{-1/2} \vec\{ \bfE(\nu) \bfX^\top(\nu)  \}
\end{align}
Slutsky's theorem and relation~(\ref{th1a}) give~(\ref{th1c}), using the following argument:
\begin{align*}
  \bfTheta(\nu)&=\left(\bfOmega^{-1}(\nu) \otimes \bfI_d\right)\sum_{h=-\infty}^{\infty}\mathbb{E}\left(\bfX_n(\nu)\bfX_{n-h}^{\top}(\nu) \otimes \bfepsilon_{ns+\nu}\bfepsilon_{(n-h)s+\nu}^{\top}\right)\left(\bfOmega^{-1}(\nu) \otimes \bfI_d\right)
\\
 &= \sum_{h=-\infty}^{\infty}\mathbb{E}\left[\bfOmega^{-1}(\nu)\bfX_n(\nu)\bfX_{n-h}^{\top}(\nu)\bfOmega^{-1}(\nu) \otimes \bfepsilon_{ns+\nu}\bfepsilon_{(n-h)s+\nu}^{\top}\right].
\end{align*}
The joint asymptotic normality of
$N^{1/2} \{ \hat{\bfbeta}^\top(1) - \bfbeta^\top(1), \ldots, \hat{\bfbeta}^\top(s) - \bfbeta^\top(s) \}$
follows using the same kind of manipulations as those for a single season $\nu$. We also hawe
\begin{eqnarray}\nonumber
\label{th1d}
  N^{1/2}\{ \hat{\bfbeta} - \bfbeta\}
  &\cd& N_{sd^2p(\nu)}\left(\bfzero, \bfTheta \right),
\end{eqnarray}
where the asymptotic covariance matrix $\bfTheta$ is a block matrix, with the asymptotic variances given by $\bfTheta(\nu)$, $\nu=1,\dots,s$, and the asymptotic covariances given by:
\[\lim_{N\to\infty}\text{cov}\left(N^{1/2}\{ \hat{\bfbeta}(\nu) - \bfbeta(\nu) \},N^{1/2}\{ \hat{\bfbeta}(\nu') - \bfbeta(\nu') \}\right)=\left(\bfOmega^{-1}(\nu) \otimes \bfI_d\right)\sum_{h=-\infty}^{\infty}\mathbb{E}\left(\bfX_n(\nu)\bfX_{n-h}^{\top}(\nu') \otimes \bfepsilon_{ns+\nu}\bfepsilon_{(n-h)s+\nu'}^{\top}\right)\left(\bfOmega^{-1}(\nu') \otimes \bfI_d\right),\]
for $\nu \neq \nu'$ and $\nu, \nu' = 1,\ldots,s$.

%\end{proof}
\subsection{Proof of Theorem \ref{estHAC}}
%\begin{proof}
Observe that
\begin{align*}
\hat{\bfPsi}^\mathrm{HAC}(\nu)-\bfPsi(\nu)&
=\sum_{h=-T_N}^{T_N}f(hb_N)\left(\hat{\Lambda}_{h}(\nu)-{\Lambda}_{h}(\nu)\right)+\sum_{h=T_N}^{T_N}\left\{f(hb_N)-1\right\}{\Lambda}_{h}(\nu)-\sum_{|h|> T_N}{\Lambda}_{h}(\nu).
\end{align*}
By the triangular inequality, for any multiplicative norm, we have
\begin{eqnarray*}
\left\|\hat{\bfPsi}^\mathrm{HAC}(\nu)-\bfPsi(\nu)\right\|&
\leq&g_1+g_2+g_3,
\end{eqnarray*} where
\begin{align*}g_1&=\sup_{|h|<N}\left\|\hat{\Lambda}_{h}(\nu)-{\Lambda}_{h}(\nu)\right\|\sum_{|h|\leq
T_N}\left|f(hb_N)\right|,&\\g_2&=\sum_{|h|\leq
T_N}\left|f(hb_N)-1\right|\left\|{\Lambda}_{h}(\nu)\right\|\quad\mbox{and}\quad g_3=\sum_{|h|>T_N}\left\|{\Lambda}_{h}(\nu)\right\|.&
\end{align*}
In view of this last inequality, to prove the convergence in probability of $\hat{\bfPsi}^\mathrm{HAC}(\nu)$ to $\bfPsi(\nu)$, it suffices to show that the probability limit of $g_1$, $g_2$ and $g_3$ is $0$.
\paragraph{$\diamond$ Step 1: convergence in probability of $\sup_{|h|<N}\left\|\hat{\Lambda}_{h}(\nu)-{\Lambda}_{h}(\nu)\right\|$ to $0$} \ \\
Let $\Lambda^{*}_h(\nu)$ be the matrix defined, for $0\leq h<N$, by
$$\Lambda^{*}_h(\nu)=\frac{1}{N}\sum_{n=0}^{N-h-1}\bfW_{n}(\nu)\bfW_{n-h}^{\top}(\nu)\quad\text{and}\quad\Lambda^{*}_{-h}(\nu)=\Lambda^{*\top}_h(\nu).$$
Observe that
$$\sup_{|h|<N}\left\|\hat{\Lambda}_{h}(\nu)-{\Lambda}_{h}(\nu)\right\|\leq \sup_{|h|<N}\left\|\hat{\Lambda}_{h}(\nu)-{\Lambda}^{*}_{h}(\nu)\right\|+\sup_{|h|<N}\left\|{\Lambda}^{*}_{h}(\nu)-{\Lambda}_{h}(\nu)\right\|.$$
By the ergodic theorem, we have
\begin{eqnarray}\label{cvLambd}
\Lambda^{*}_h(\nu)&\cas&\Lambda_h(\nu).
\end{eqnarray}
A Taylor expansion of $\vec\{\hat{\Lambda}_h(\nu)\}$ around $\bfbeta$
and \eqref{th1a} give
\begin{eqnarray*}
\vec\{\hat{\Lambda}_h(\nu)\}=\vec\{\Lambda^{*}_h(\nu)\}+\frac{\partial \vec\{\Lambda^{*}_h(\nu)\}}{\partial\bfbeta^\top(\nu)}
(\hat\bfbeta(\nu)-\bfbeta(\nu))+\mathrm{O}_\mathbb{P}\left(\frac{1}{N}\right).
\end{eqnarray*}
In view of \eqref{cvLambd} and by {\bf (A3)}, we then deduce that
\begin{equation}\label{dergammaborn}
\lim_{N\to\infty}\sup_{|h|<N}\left\|\frac{\partial \vec\{\Lambda^{*}_h(\nu)\}}{\partial\bfbeta^\top(\nu)}\right\|<\infty,\quad a.s.
\end{equation}
%$$b_n\sum_{|h|<N}|f(hb_N)|=O(1).$$
By the ergodic theorem, \eqref{th1a} and \eqref{dergammaborn}, for any multiplicative norm, we have
\begin{eqnarray}
\sup_{|h|<N}\left\|\vec\left(\hat{\Lambda}_h(\nu)-\Lambda^{*}_h(\nu)\right)\right\|\leq\lim_{N\to\infty}\sup_{|h|<N}\left\|\frac{\partial \vec\{\Lambda^{*}_h(\nu)\}}{\partial\bfbeta^\top(\nu)}\right\| \left\|
\hat\bfbeta(\nu)-\bfbeta(\nu)\right\|+\mathrm{O}_\mathbb{P}\left(\frac{1}{N}\right)=\mathrm{O}_\mathbb{P}\left(\frac{1}{\sqrt{N}}\right).\label{cvLambd2}
\end{eqnarray}
From \eqref{cvLambd} and \eqref{cvLambd2}, we deduce that
\begin{equation}
\label{cvLambda}
\sup_{|h|<N}\left\|\hat{\Lambda}_{h}(\nu)-{\Lambda}_{h}(\nu)\right\|=\mathrm{O}_\mathbb{P}\left(\frac{1}{\sqrt{N}}\right)=\mathrm{o}_{\mathbb{P}}(1),
\end{equation}
the conclusion is complete.

\paragraph{$\diamond$ Step 2: convergence in probability of $g_1$, $g_2$ and $g_3$ to $0$} \ \\
By {\bf (A3)}, $\mathbb{E}\|\bfW_{n}\|^{2+\kappa}<\infty$. Davydov's inequality (see \cite{D68}) then entails that
\begin{equation}
\label{M_ijhborn}
\left\|{\Lambda}_{h}(\nu)\right\|=\left\|\cov\left(\bfW_{n}(\nu),\bfW_{n-h}(\nu)\right)\right\|\leq K\alpha_{\bfepsilon}^{\kappa/(2+\kappa)}([h/2]).
\end{equation}
 In view of {\bf (A3)}, we thus have $g_3\to 0$ as
$N\to\infty$. Let $m$ be a fixed integer and we write $g_2\leq
s_1+s_2$, where $$s_1=\sum_{|h|\leq
m}\left|f(hb_N)-1\right|\left\|{\Lambda}_{h}(\nu)\right\|\quad\text{and}\quad s_2=\sum_{m<|h|\leq
T_N}\left|f(hb_N)-1\right|\left\|{\Lambda}_{h}(\nu)\right\|.$$ For
$|h|\leq m$, we have $hb_N\to0$ as $N\to\infty$ and $f(hb_N)\to1$,
it follows that $s_1\to0$. If we choose $m$ sufficiently large,
$s_2$ becomes small. Using \eqref{M_ijhborn} and the fact that
$f(\cdot)$ is bounded, it follows that $g_2\to0$.

In view of \eqref{fenetreO(1)} and \eqref{cvLambda}, we have
\begin{align*}g_1&=\sup_{|h|<N}\left\|\hat{\Lambda}_{h}(\nu)-{\Lambda}_{h}(\nu)\right\|\sum_{|h|\leq
T_N}\left|f(hb_N)\right|,\\ &= \frac{1}{b_N}\sup_{|h|<N}\left\|\hat{\Lambda}_{h}(\nu)-{\Lambda}_{h}(\nu)\right\|b_N\sum_{|h|\leq
T_N}\left|f(hb_N)\right|,
\\ &\leq \frac{1}{b_N}\sup_{|h|<N}\left\|\hat{\Lambda}_{h}(\nu)-{\Lambda}_{h}(\nu)\right\|\mathrm{O}(1)=\mathrm{O}_\mathbb{P}\left(\frac{1}{b_N\sqrt{N}}\right)=\mathrm{o}_\mathbb{P}(1),
\end{align*}
since $Nb_N^2\to\infty$, in view of (\ref{conditions_bn}). The proof is complete.
%\end{proof}

\noindent
{\bf Acknowledgements:}
We sincerely thank the
anonymous reviewers and editor for helpful remarks.

%% Authors are advised to submit their bibtex database files. They are
%% requested to list a bibtex style file in the manuscript if they do
%%% not want to use model1a-num-names.bst.

%%%%% Debut Bibliographie
%\bibliographystyle{apalike}
%\bibliography{weak_PVAR}
%\addcontentsline{toc}{section}{References}
%%%%% Fin Bibliographie

\end{document}